\documentclass[journal,twoside,web]{ieeecolor}
\usepackage{generic}
\usepackage{epsfig} 

\usepackage{cite}
\usepackage{amsmath,amssymb,amsfonts}

\makeatletter
\let\proof\@undefined                        
\let\endproof\@undefined                  
\makeatother
\usepackage{amsthm}

\usepackage{algorithmic}
\usepackage{graphicx}
\usepackage{transparent} 
\usepackage{algorithm,algorithmic}
\usepackage{dsfont} 
\usepackage{bm} 
\usepackage{makecell} 
\usepackage{tikz} 
\usepackage{subcaption}

\allowdisplaybreaks

\usepackage{hyperref}
\hypersetup{
    colorlinks,
    linkcolor={blue},
    citecolor={blue},
    urlcolor={blue}
} 
\usepackage{textcomp}
\def\BibTeX{{\rm B\kern-.05em{\sc i\kern-.025em b}\kern-.08em
    T\kern-.1667em\lower.7ex\hbox{E}\kern-.125emX}}
\theoremstyle{plain} 
\newtheorem{definition}{\textbf{Definition}}
\newtheorem{assumption}{\textbf{Assumption}}
\newtheorem{theorem}{\textbf{Theorem}}
\newtheorem{lemma}{\textbf{Lemma}}
\newtheorem{proposition}{\textbf{Proposition}}
\newtheorem{corollary}{\textbf{Corollary}}
\theoremstyle{definition} 
\newtheorem{remark}{\textbf{Remark}}

\renewcommand{\QED}{\QEDopen}

\newcommand{\qt}[1]{{\texttransparent{1}{\color{black}#1}}}

\newcommand{\bR}{\mathbb{R}}
\newcommand{\bZ}{\mathbb{Z}}
\newcommand{\bN}{\mathbb{N}}
\newcommand{\bI}{\mathbb{I}}

\newcommand{\cX}{\mathcal{X}}
\newcommand{\cU}{\mathcal{U}}

\newcommand{\cS}{\mathcal{S}}

\newcommand{\cK}{\mathcal{K}}

\newcommand{\xp}[1]{x^+(#1)}

\newcommand{\splus}{\hspace{-0.08cm}+\hspace{-0.08cm}}
\newcommand{\sminus}{\hspace{-0.06cm}-\hspace{-0.06cm}}

\newcommand{\sleq}{\hspace{-0.08cm}\leq\hspace{-0.08cm}}

\newcommand{\bzero}{\mathbf{0}}
\newcommand{\bone}{\mathbf{1}}
\newcommand{\bfI}{\mathbf{I}}
\newcommand{\bigo}[1]{\mathcal{O}\left(#1\right)}

\newcommand{\compratio}[1]{\mathcal{R}_{\mathrm{cr},#1}}
\newcommand{\compratiolq}[1]{\mathcal{R}_{\mathrm{cr},\mathrm{LQ},#1}}

\newcommand{\para}{\theta}
\newcommand{\mpara}{\delta(\theta)}
\newcommand{\epara}{\varepsilon_{\theta}}
\newcommand{\tpara}{\theta^\ast}
\newcommand{\hpara}{\hat{\theta}}
\newcommand{\setpara}{\Theta(\hat{\theta}, \varepsilon_{\theta})}

\newcommand{\dsys}{\Delta_f(x,u;\theta)}

\newcommand{\cost}{\ell}

\newcommand{\probopt}[1]{\mathrm{P}_{\mathrm{IHOPC}}(x;{#1})}
\newcommand{\probmpc}[1]{\mathrm{P}_{\mathrm{MPC}}(x;{#1})}

\newcommand{\lyacf}{\varepsilon_F}
\newcommand{\umpc}{\overline{\gamma}}
\newcommand{\mpc}{\gamma}

\newcommand{\ussbh}{\overline{\sigma}_{\mathbf{H}}}
\newcommand{\ussbr}{\overline{\sigma}_{\mathbf{R}}}
\newcommand{\lssbh}{\underline{\sigma}_{\mathbf{H}}}
\newcommand{\edsc}{H}
\newcommand{\edsr}{\lambda}

\newcommand{\nA}{\nu_A}
\newcommand{\nB}{\nu_B}

\newcommand{\eab}{\varepsilon_{A,B}}

\newcommand{\bu}{\mathbf{u}}

\newcommand{\otx}[1]{\psi_{x, #1}}

\newcommand{\ctx}[1]{\phi_{x, #1}}

\newcommand{\ctu}[1]{\phi_{u, #1}}

\newcommand{\sueds}[1]{\omega_{#1,N}}
\newcommand{\leds}[1]{\eta_{#1,N}}

\newcommand{\xroa}{\mathcal{X}_{\mathrm{ROA}}}

\newcommand{\nQ}[1]{\| #1 \|^2_Q}
\newcommand{\nR}[1]{\| #1 \|^2_R}
\newcommand{\emax}[1]{\lambda_{\max}(#1)}
\newcommand{\emin}[1]{\lambda_{\min}(#1)}

\newcommand{\ldd}{L_{\Delta}}
\newcommand{\lddx}{L_{\Delta,x}}
\newcommand{\lddu}{L_{\Delta,u}}

\newcommand{\difxi}[1]{\|\Delta \xi^\star_{#1}(x)\|}

\newcommand{\difu}[1]{\|\Delta u^\star_{#1}(x)\|}
	
\newcommand{\dpairf}[1]{\|\Delta^\star_{f,#1}(x;\hpara)\|}

\newcommand{\dpaf}[1]{\|\Delta^\star_{f,#1}(x;\para)\|}

\newcommand{\dxip}[1]{\|\Delta \xi^+_{#1}\|}
\newcommand{\dup}[1]{\|\Delta u^+_{#1}\|}
\newcommand{\xiph}[1]{\xi^\star_{#1}(x^+(\hpara);\hpara)}
\newcommand{\uph}[1]{u^\star_{#1}(x^+(\hpara);\hpara)}
\newcommand{\law}{\mu_{N,\hpara}(x)}

\newcommand{\xin}[1]{\xi^\star_{#1}(x;\hpara)}
\newcommand{\xit}[1]{\xi^\star_{#1}(x;\para)}
\newcommand{\uun}[1]{u^\star_{#1}(x;\hpara)}
\newcommand{\uut}[1]{u^\star_{#1}(x;\para)}

\newcommand{\omelq}{\overline{\omega}^\ast_{\mathrm{LQ}}}

\newcommand{\mec}{M_{\mathrm{ec}}}
\newcommand{\blfx}{\overline{L}_{f,x}}
\newcommand{\blfu}{\overline{L}_{f,u}}
\newcommand{\llx}{L_{\ell,x}}
\newcommand{\llu}{L_{\ell,u}}
\newcommand{\mlx}{m_{\ell,x}}
\newcommand{\mlu}{m_{\ell,u}}
\newcommand{\bo}{\overline{\omega}}

\begin{document}
\title{Certainty-Equivalence Model Predictive Control: Stability, Performance, and Beyond}
\author{Changrui Liu, \IEEEmembership{Graduate Student Member, IEEE}, Shengling Shi, and Bart De Schutter, \IEEEmembership{Fellow, IEEE}
\thanks{This paper is part of a project that has received funding from the European Research Council (ERC) under the European Union’s Horizon 
2020 research and innovation programme (Grant agreement No. 101018826 - CLariNet). (Corresponding authors: Changrui Liu \& Shengling Shi)}
\thanks{Changrui Liu and Bart De Schutter are with the Delft Center for Systems and Control, TU Delft, Delft, The Netherlands. (e-mail: \{c.liu-14, b.deschutter\}@tudelft.nl). }
\thanks{Shengling Shi was with the Department of Chemical Engineering, MIT, Cambridge, United States, and is now with Delft Center for Systems and Control, TU Delft. (e-mail: shengling.shi@tudelft.nl).}
}

\maketitle

\begin{abstract}
    Handling model mismatch is a common challenge in model predictive control (MPC). While robust MPC is effective, its conservatism often makes it less desirable. Certainty-equivalence MPC (CE-MPC), which uses a nominal model, offers an appealing alternative due to its design simplicity and low computational costs. This paper investigates CE-MPC for uncertain nonlinear systems with multiplicative parametric uncertainty and input constraints that are inactive at the steady state. The primary contributions are two-fold. First, a novel perturbation analysis of the MPC value function is provided, without assuming the Lipschitz continuity of the stage cost, better tailoring the widely used quadratic cost and having broader applicability in value function approximation, learning-based MPC, and performance-driven MPC design. Second, the stability and performance analysis of CE-MPC are provided, quantifying the suboptimality of CE-MPC compared to the infinite-horizon optimal controller with perfect model knowledge. The results provide insights in how the prediction horizon and model mismatch jointly affect stability and the worst-case performance. Furthermore, the general results are specialized to linear quadratic control, and a competitive ratio bound is derived, serving as the first competitive-ratio bound for MPC of uncertain linear systems with input constraints and multiplicative uncertainty.
\end{abstract}

\begin{IEEEkeywords}
Predictive control, optimal control, uncertain systems, performance analysis, perturbation analysis
\end{IEEEkeywords}
\vspace*{-0.5cm}
\section{Introduction}
\label{sec:introduction}
\IEEEPARstart{M}{odel} predictive control (MPC) has prospered in many engineering fields, such as chemical processes \cite{larsen2006industrial}, robotics \cite{worthmann2015model}, and microgrids \cite{hu2021model}. In general, MPC requires a model to predict future behaviors of the system to computie the optimal inputs. Therefore, the effectiveness of MPC significantly relies on the model accuracy \cite[Section 1.2]{rawlings2017model}, and successful implementations usually require model tuning \cite{qin2003survey}. Nonetheless, obtaining a perfect model is impossible due to inherent limitations in modeling and identification (e.g., unmodeled dynamics, limited data, and measurement noise), and hence, model mismatch always exists in practice.

For MPC, model mismatch is commonly described using a nominal model with \qt{general} additive disturbances \cite{cannon2011robust, mayne2011tube, manzano2020robust}. To quantify parametric uncertainty, additive residual models have been frequently adopted \cite{liu2025robust, adetola2009adaptive, kohler2020computationally, dogan2023regret}. 
Another popular method considers parametric uncertainty within a structured nominal model, e.g., parametric linear models \cite{fleming2014robust, moreno2022performance}. There have been recent trends on studying the joint effect of multiplicative parametric uncertainty and additive disturbances for linear systems \cite{calafiore2012robust, lorenzen2019robust, arcari2023stochastic, shin2022near}. Besides, results on parametric nonlinear systems \cite{wabersich2022cautious, lin2022bounded} also exist. To handle model mismatch, most contributions are about robust MPC \cite{mayne2011tube, schwenkel2020robust, kohler2020computationally, lorenzen2019robust} or MPC with model learning \cite{dogan2023regret, moreno2022performance, wabersich2022cautious}, focusing on stability and feasibility, whereas only a few studies provide performance guarantees in terms of (sub)optimality \cite{li2023certifying, dogan2023regret, moreno2022performance, schwenkel2020robust}. Moreover, results on uncertain nonlinear systems without using a linear nominal model are scarce \cite{schwenkel2020robust, lin2022bounded}.

Optimality-related guarantees for MPC have been studied in terms of the \textit{infinite-horizon} performance using relaxed dynamic programming (RDP) \cite{grune2008infinite} (see also a compact discussion in \cite[Section 2.2.3]{mayne2014model}). Although performance analysis using RDP is quite mature, most of the existing results are only valid when the model is perfect \cite{grune2010analysis, kohler2021stability, kohler2023stability}. Some recent extensions have covered nonlinear systems with disturbances \cite{schwenkel2020robust}, a class of parameterized linear systems with polytopic uncertainty \cite{moreno2022performance}, and general uncertain linear systems \cite{liu2024stability, shi2025suboptimality}. In \cite{schwenkel2020robust}, strong assumptions on the existence of some unjustified growth functions have been made to establish the performance bound \cite[Proposition 2]{schwenkel2020robust}, and the result only applies to Lipschitz-continuous costs. Besides, the analysis in \cite{moreno2022performance} is not consistent, i.e., the obtained bound does not degenerate to the standard results (see, e.g., \cite{grune2017nonlinear}) when no model mismatch exists. The performance bound in \cite{liu2024stability, shi2025suboptimality} preserves consistency and it does not require strong assumptions, but it is still tailored for linear systems. Therefore, there is still a lack of meaningful results on infinite-horizon performance analysis of MPC for uncertain nonlinear systems.

Suboptimality due to model mismatch has also been quantified using \textit{regret} \cite{dogan2023regret, lin2022bounded} and \textit{competitive ratio} \cite{li2023certifying, goel2022competitive}. However, these results, being finite-horizon in nature, cannot reveal how MPC approximates infinite-horizon optimal control, thus advocating for RDP-based methods. Furthermore, obtaining perfect models is impossible even with model learning \cite{dogan2023regret, wabersich2022cautious, moreno2022performance, lorenzen2017adaptive}. Therefore, it is still fundamental and necessary to study MPC in the certainty-equivalence framework \cite{hespanha1999certainty}, known as certainty-equivalence MPC (CE-MPC) \cite{soloperto2019dual}. This work focuses on \textit{multiplicative} parametric uncertainty without disturbances or state estimation error, which is different from the settings in \cite{grimm2007nominally, allan2017inherent, pannocchia2011conditions}, where additive disturbances and state estimation errors are considered. Besides, we only consider the common scenarios where input constraints are \textit{inactive} at the steady state \cite{mayne2000constrained}.

The stability and performance of CE-MPC for uncertain nonlinear systems remains largely unexplored. The joint effect of the prediction horizon and the model mismatch has only been studied for finite-state Markov decision processes with bounded rewards \cite{zhang2024predictive}. In terms of stability, CE-MPC has recently been investigated for parametric nonlinear systems in \cite{kuntz2024beyond}, which, though being related to the stability results of this paper, is based on a different setting, where the true model can be time-varying, terminal ingredients are present, and the parameters of the nominal model are \qt{set to} \textit{zero}. Besides, in \cite{kuntz2024beyond}, the existence of \qt{some $\mathcal{K}$ functions (e.g., \cite[Assumption 9]{kuntz2024beyond}) to quantify the influence of parametric mismatch is assumed to prove stability}, which is not needed in the current work. In summary, we will address the following questions:
\begin{enumerate}
	\item \textit{Under what conditions is the CE-MPC controller designed based on the nominal model guaranteed to stabilize the true model at the true equilibrium?}
	\item \textit{In terms of the infinite-horizon performance, how suboptimal is the CE-MPC controller compared to the offline infinite-horizon controller that knows the true model?}
	\item \textit{What is the interplay between the prediction horizon and the model mismatch in terms of stability and (worst-case) performance bounds of CE-MPC?}
\end{enumerate}
In response to the above questions, this work advances the state of the art through the following contributions:
\begin{enumerate}
    \item We present an analysis pipeline for CE-MPC of nonlinear systems with multiplicative parametric uncertainty and input constraints, resulting in a novel suboptimality bound to quantify its infinite-horizon performance with respect to the (offline) oracle infinite-horizon optimal controller (IHOPC) that knows the true model.
    \item A novel perturbation analysis of the MPC value function is provided without using the Lipschitz continuity of the cost, which can be of independent interest for learning-based MPC. The perturbation analysis helps establish the suboptimality bounds that reveal the joint effect of the prediction horizon and the model mismatch on the stability and the worst-case infinite-horizon performance. 
    \item Denote the true parameter and the nominal parameter by $\tpara$ and $\hpara$, respectively. An upper bound of the normed model mismatch $\|\tpara - \hpara\|$ is derived to ensure closed-loop stability. Besides, competitive-ratio upper bound $\mathcal{R}$ is obtained under mild conditions, i.e., $\mathcal{R}$ satisfies
    \begin{equation*}
    	\label{eq:competitive_introduction}
    	\mathcal{R} \geq \sup_{\|\tpara-\hpara\|\leq \epara}\frac{\text{Cost of CE-MPC}(\hpara)}{\text{Cost of IHOPC}(\tpara)},
    \end{equation*}
    where $\epara > 0$ is the considered worst-case mismatch level. Moreover, $\mathcal{R}$ explicitly depends on the prediction horizon of CE-MPC, and an \textit{approximate} optimal horizon can be computed by minimizing $\mathcal{R}$. This horizon, though not necessarily optimal for true performance, can be used as a coarse warm start for horizon tuning when lacking perfect system knowledge.
    \item We further specify the performance bounds for uncertain linear systems with quadratic costs, extending our previous results \cite{liu2024stability} such that a \textit{pure-ratio} bound is obtained.
\end{enumerate}

Importantly, our analysis is consistent with those without model mismatch \cite{grune2008infinite, grune2017nonlinear, kohler2023stability}, since our results still advocate for a longer prediction horizon when no mismatch exists. It should be emphasized that we do not develop new MPC strategies, but focus on theoretical analysis of CE-MPC. The rest of the paper is organized as follows. Section \ref{sec:preliminaries} provides preliminaries. The problems of analyzing the stability and the infinite-horizon performance of CE-MPC are formulated in Section \ref{sec:problem_formulation}. Perturbation analysis of MPC value function is given in Section \ref{sec:theoretical_analysis}, and Section \ref{sec:theoreticalAnalysis:stability_suboptimality} presents our main results. Simulations and extension to linear systems are presented in Section \ref{sec:examples}. Finally, we conclude the paper in Section \ref{sec:conclusion}.

\section{Preliminaries}
\label{sec:preliminaries}
\subsection{Notation}
$\bR$ denotes the set of real numbers, and $\bR_+$ is the set of non-negative real numbers. The set of natural numbers is $\bN$, the set of integers is $\bZ$, and $\bI_{a:b} := \bN \cap [a, b]$ for $0 \leq a \leq b$. For a vector $x$, its $i$-th element is $x[i]$, and the sequence $\{x_k\}^{t_2}_{k= t_1}$ is denoted by $\mathbf{x}_{t_1:t_2}$. The bold symbols $\bzero$, $\bone$, and $\bfI$ denote the zero vector, one vector, and identity matrix of appropriate dimensions. The ceiling function\footnote{For $a \in \bR$, $\lceil \cdot \rceil$ is defined as $\lceil a \rceil := \min\{n \in \bZ \mid n \geq a\}$.} is denoted by $\lceil\cdot\rceil$. 

Let $M$ be a matrix. Then its transpose is $M^\top$, and $M \geq (>) \;0$ and $M \succ (\succeq) \;0$ denote element-wise non-negativity (positivity) and positive (semi-) definiteness, respectively. $\|\cdot\|$ denotes the $2$-norm of a vector and the induced $2$-norm of a matrix. The Euclidean ball centered at $x$ with radius $R$ is denoted by $\mathcal{B}(x, R) := \{x' \in \bR^{n}\mid \|x - x'\| \leq R\}$. For a \textit{symmetric} matrix $M$, its largest and smallest eigenvalues are denoted, respectively, by $\emax{M}$ and $\emin{M}$. Moreover, $r(M) := \emax{M}/\emin{M}$, and $\|x\|_{M}$ stands for $(x^\top M x)^{1/2}$. The $i$-th row of $M$ is denoted by $[M]_{(i,:)}$. 

Given a set $\cS$, $\cS^{n+1} := \cS^{n} \times \cS$ for $n \in \bI_{1:\infty}$ with $\cS^1 := \cS$, where $\times$ is the Cartesian product. A function $\alpha: [0, a) \to \bR_+$ belongs to class $\mathcal{K}$ (i.e., $\alpha \in \cK$) if it is strictly increasing and $\alpha(0) = 0$. $\nabla f(x)$ is the gradient of $f$ at $x$. A function $f: \bR^n \to \bR_+$ is positive definite if $f(x) = 0 \iff x = \bzero$. 
\subsection{Infinite-horizon Optimal Control}
Consider discrete-time nonlinear systems modeled by
\begin{equation}
    \label{eq:sec_pre:nonlinear_system}
    x_{t+1} = f(x_t, u_t; \para),
\end{equation}
where $t \in \bN$ is the time step, $x_t \in \cX = \bR^n$ and $u_t \in \cU \subseteq \bR^m$ are, respectively, the state and input at time step $t$, $f: \cX \times \cU \to \cX$ is the model function, and $\para \in \bR^{n_{\theta}}$ is the model parameter. In addition, the set $\cU$ is \textit{closed} and is defined by $g_u: \bR^{m} \to \bR^{c_u}$ as $\cU = \{u \in \bR^m \mid g_u(u) \leq \mathbf{1}\}$. The \textit{true} model is defined by \eqref{eq:sec_pre:nonlinear_system} with $\para = \tpara$, whereas a practical controller only has access to nominal parameter $\hpara$ due to lack of precise knowledge or estimation errors, which are prevalent in identification and model learning \cite{arcari2023stochastic, lin2022bounded}. The metric
\begin{equation}
	\label{eq:para_metric}
	\mpara = \|\hpara - \para\|
\end{equation}
is used to quantify the difference between $\para$ and $\hpara$. Although $\tpara$ is \textit{unknown}, we assume that $\tpara$ belongs to a known compact set $\setpara \subset \bR^{n_{\para}}$ defined by
\begin{equation}
	\label{eq:sec_pre:parameter_set}
	\setpara := \{\para \in \bR^{n_{\para}} \mid \|\hpara - \para\| \leq \epara\},
\end{equation}
where $\epara$ is the \textit{mismatch level}. 
\begin{assumption}
	\label{ass:sec_pre:basic_problem}
	$f$ and $g_u$ satisfy the following conditions:
	\begin{enumerate}
		\item $f$ and $g_u$ are twice continuously differentiable, and for all $(x,u) \in \cX \times \cU$, the function $f_{x,u}: \bR^{n_\theta} \to \cX$ defined by $f_{x,u}(\theta) := f(x, u;\theta)$ is also twice continuously differentiable;
		\item $f(\bzero, \bzero;\theta) = \bzero$ for all $\para \in \setpara$;
		\item $g_u(\bzero) \leq \mathbf{1}$ (i.e., $u = \bzero$ lies in the interior of $\cU$).
	\end{enumerate}
\end{assumption}
The first condition is standard for sensitivity analysis in control \cite{lin2022bounded} and optimization \cite{shin2022exponential}. The second condition states that model mismatch does not change the true equilibrium $(x,u)=(\bzero,\bzero)$, such that the CE-MPC controller is well informed of the true equilibrium. Kernel-based models \cite{van2022kernel} suit this condition, and so do many applications (e.g., \qt{some} HVAC systems \cite{dogan2023regret}, cruise control \cite{ames2016control}, and robot arms \cite{li2025learning}). However, systems with additive disturbances or bias violate this condition. The last condition ensures that the input constraints are \textit{inactive} at $(x,u) = (\bzero, \bzero)$. For convenience, we denote $f(x, u; \para)$ occasionally by $x^+(u;\para)$, and for any $\para \in \setpara$ and $(x, u) \in \cX \times \cU$, we define the model error $\dsys$ as
\begin{equation}
	\label{eq:dynamic_parametric_difference}
	\dsys := f(x, u;\para) - f(x, u;\hpara),
\end{equation}
and it holds that $\Delta_f(\mathbf{0},\mathbf{0};\para) = \mathbf{0}$ due to Assumption 1, item 2). In addition, we define
\begin{equation}
	R(x;\epara) := \max_{u \in \mathcal{U}, \para \in \setpara}\|\dsys\|,
\end{equation}
which is the maximum one-step-ahead state deviation due to model mismatch given the current state $x$. At each time step $t$, the system under control incurs a non-trivial\footnote{In this work, trivial constant cost functions $\ell \equiv c \;(c \geq 0)$ are excluded since they have no effect on regulation or stabilization.} positive-definite stage cost $\cost(x_t, u_t) \geq 0$ quantified by $\cost: \bR^n \times \bR^m \to \bR_+$, and the optimized stage cost $\cost^\star: \bR^n \to \bR_+$ is defined as
\begin{equation}
	\label{eq:sec_pre:optimized_stage_cost}
	\cost^\star(x) := \min_{u \in \cU}\cost(x, u).
\end{equation}
which will be useful for performance analysis.

For the model in \eqref{eq:sec_pre:nonlinear_system}, starting from $x = x_t$, the $k$-step-ahead \textit{open-loop} state under the sequence $\bu_{t:t+K-1} = \{u_{t+i}\}^{K-1}_{i=0}$ is denoted by $\otx{\para}(k;x,\bu_{t:t+K-1})$ for all $k \in \bI_{0:K}$, where $\otx{\para}(0;x,\bu_{t:t+K-1}) = x$. Given a control policy $\pi: \cX \to \cU$, the \textit{closed-loop} state at time step $t + k$, starting from $x = x(t)$, is denoted by $\ctx{\para}(k; x,\pi)$ with $\ctx{\para}(0;x,\pi) = x$ and $\ctx{\para}(k+1;x,\pi) = f(\ctx{\para}(k;x,\pi), \pi(\ctx{\para}(k;x,\pi)); \para)$, where the associated control input $\pi\left(\ctx{\para}(k;x,\pi)\right)$ is further denoted by $\ctu{\para}(k;x,\pi)$ for all $k \in \bN$. 

The goal is to determine $\bu_{0:\infty} \in \cU^{\infty}$ that asymptotically stabilizes the system at $(x, u) = (\bzero, \bzero)$ while minimizing the \textit{infinite-horizon} cost $J_{\infty}(x_0; \bu_{0:\infty}, \para) = \sum^{\infty}_{t = 0}\cost\left(x_t, u_t\right)$ subject to the model in \eqref{eq:sec_pre:nonlinear_system}. Equivalently, starting from $x_0 = x$, the objective is to solve the following  optimization problem:
\begin{align*} 
	\probopt{\para}: & \min_{\{u_{k|0}, x_{k|0}\}^{\infty}_{k=0}} \sum^{\infty}_{k=0} \cost\left(x_{k|0}, u_{k|0}\right)
	\\ \text{s.t. } & x_{k+1|0} = f(x_{k|0}, u_{k|0}; \para), \quad \forall k \in \bN,
	\\ & \; g_u(u_{k|0}) \leq \mathbf{1}, \quad \forall k \in \bN ,
	\\ & \; x_{0|0} = x,
\end{align*}
where $x_{k|0}$ and $u_{k|0}$ are the predicted $k$-step-ahead state and input starting from time step $0$, respectively. We assume $\probopt{\para}$ admits a minimizer, denoted by $\bu^\star_{0:\infty}(x, \para)$. Moreover, the value $J_{\infty}(x; \bu^\star_{0:\infty}(x, \para), \para)$ is denoted by $V_{\infty}(x;\para)$. The controller defined by $\probopt{\tpara}$ is called the \textit{oracle} controller, the existence of which requires defining the region of attraction (ROA) as follows:
\begin{definition}[Region of Attraction]
	\label{def:sec_pre:roa}
	The set of states that can be regulated to $x = \bzero$ asymptotically is defined as $\xroa = \{x \in \cX \mid \exists \bu_{0:\infty} \in \cU^{\infty},  \mathrm{s.t.}\;\; J_{\infty}(x;\bu_{0:\infty},\tpara) < +\infty\}$.
\end{definition}
Definition \ref{def:sec_pre:roa} implicitly requires the existence of an input sequence $\bu_{0:\infty}$ such that the infinite-time cost is \textit{summable} when $(x_t,u_t)$ approaching $(\bzero,\bzero)$. Following the definition of ROA, a standard result follows immediately as in \cite{rawlings2017model}:
\begin{corollary}
	\label{corollary:sec_pre:invariant_roa}
	Under Assumption \ref{ass:pre:initial_state_roa}, $\xroa$ is a forward control invariant set, i.e., for all $x \in \xroa$, there exists a $u \in \cU$ such that $f(x,u; \tpara) \in \xroa$.
\end{corollary}
Corollary \ref{corollary:sec_pre:invariant_roa} ensures $x_t \in \xroa$ for all $t \in \bN$ under a proper control policy without imposing state constraints. Next, we present an additional assumption to ensure that the optimal control problem is well-defined such that regulating to the equilibrium $(x, u) = (\bzero, \bzero)$ is achievable.
\begin{assumption}
	\label{ass:pre:initial_state_roa}
	The initial state $x_0$ lies in $\xroa$.
\end{assumption}
Assumption \ref{ass:pre:initial_state_roa} ensures that $x_0$ can be stabilized, which is necessary since stabilizing the system may not be possible for an arbitrary $x_0$ under input constraints. Consequently, the state can always to kept within $\xroa$, albeit we do not consider explicit state constraints. Similar assumptions have been made in \cite{kohler2021stability} and \cite{kohler2023stability}, assuming the true model is known. However, whenever explicit state constraints are present, Definition \ref{def:sec_pre:roa} should be modified to include them explicitly. In the literature, if the true model is given, $\xroa$ can be characterized using sublevel sets \cite{kohler2021stability, boccia2014stability, darup2015missing} (see also the discussions in \cite[Remark 1]{kohler2023stability}). Characterizing $\xroa$ without knowing the true model, however, is out of the scope of this paper.

\section{Problem Formulation}
\label{sec:problem_formulation}
Solving $\probopt{\para}$ is \textit{intractable} in general, and MPC is an approximation to $\probopt{\para}$ \cite{grune2017nonlinear}. The current paper considers MPC without terminal ingredients, which is common in theory \cite{schwenkel2020robust, grune2017nonlinear} and practice \cite{worthmann2015model, hu2021model}. At time step $t$, given $x_t = x$, the objective is to minimize $J_{N}(x;\bu_{0|t:N-1|t}, \para) = \sum^{N-1}_{k=0}\cost(x_{k|t}, u_{k|t}) + \ell^\star(x_{N|t})$ by solving
\begin{align*} 
\probmpc{\para}: & \min_{\{u_{k|t}\}^{N-1}_{k=0}, \{x_{k|t}\}^{N}_{k=0}} \sum^{N-1}_{k=0} \cost(x_{k|t}, u_{k|t}) + \ell^\star(x_{N|t}) 
\\ \text{s.t. } & x_{k+1|t} = f(x_{k|t}, u_{k|t}; \para), \quad \forall k \in \mathbb{I}_{0:N-1}, 
\\ & \; g_u(u_{k|t}) \leq \mathbf{1}, \quad \forall k \in \mathbb{I}_{0:N-1}, 
\\ & \; x_{0|t} = x_t, 
\end{align*}
where $N \in \bI_{1:\infty}$ is the prediction horizon, $x_{k|t}$ and $u_{k|t}$ are, respectively, the $k$-step-ahead predicted state and input from time step $t$. In this work, we consider a class of cost functions that are commonly used in optimal and predictive control.  
\begin{assumption}[Well-conditioned Costs]
	\label{ass:sec_prob:terminal_cost}
	The function $\ell$ is separable as $\ell(x,u) = \ell_x(x) + \ell_u(u)$, and there exist scalars $L_{\ell,x}$, $L_{\ell,u}$, $m_{\ell,x}$, and $m_{\ell,u}$ such that $\ell_x$ and $\ell_u$ satisfy the following conditions:
	\begin{enumerate}
		\item $\ell_x$ and $\ell_u$ are twice continuously differentiable;
		\item $\ell_x$ is $L_{\ell,x}$-strongly smooth\footnote{A function $f: \mathcal{D} \to \mathbb{R}$ is $L$-strongly smooth if for any $x',x'' \in \text{Int}(\mathcal{D})$, it holds that $\|\nabla f(x') - \nabla f(x'')\| \leq L\|x'-x''\|$.} and $m_{\ell,x}$-strongly convex\footnote{A function $f: \mathcal{D} \to \mathbb{R}$ is $m$-strongly convex if for any $x', x'' \in \mathcal{D}$, it holds that $f(x'')\geq f(x') + \nabla f(x')^\top (x''-x') + \frac{m}{2}\|x'-x''\|$.};
		\item $\ell_u$ is $L_{\ell,u}$-strongly smooth and $m_{\ell,u}$-strongly convex;
		\item $\ell_x$ and $\ell_u$ are non-negative, $\ell_x(\mathbf{0}) = \ell_u(\mathbf{0}) = 0$, and $\nabla\ell_x(\mathbf{0}) = \nabla\ell_x(\mathbf{0}) = \mathbf{0}$;
		\item $\ell_x$ is a relaxed control Lyapunov function (CLF), i.e., there exists a constant $\lyacf \in \bR_{+}$ such that
		\begin{equation}
			\label{eq:sec_prob:relaxed_clf}
			\min_{u \in \cU}\{\ell_x(x^+(u;\hpara)) + \cost(x, u)\} \leq (1 + \lyacf)\ell_x(x).
		\end{equation}
	\end{enumerate}
\end{assumption}
\begin{corollary}
	\label{coro:optimized_stage_cost}
	Let Assumption \ref{ass:sec_prob:terminal_cost} hold. Then $\ell^\star = \ell_x$.
\end{corollary}
First, the separable cost condition is common in both theory \cite{lin2021perturbation} and practice \cite{dall2016optimal}, which, however, rules out rate of change penalties on the input. The differentiability condition is essential for perturbation analysis  \cite{bonnans2013perturbation} since it ensures a differentiable objective of $\probmpc{\para}$, guaranteeing the existence of a minimizer given that $\cU^N$ is compact~\cite{rawlings2017model}. The smoothness and convexity conditions are common in optimal control \cite{lin2021perturbation, li2021online, li2020online2, shi2020online} and convex optimization \cite{nesterov2013introductory}. Condition 4) is less common, but can be satisfied via re-parameterization \cite{lin2021perturbation} and power lifting (i.e., consider $\ell'(x) = \ell^2(x)$, knowing $\ell(\mathbf{0}) = 0$) without loss of generality. Finally, the relaxed CLF condition 5) is used to study the stability of the nominal system \cite{kohler2023stability, kohler2021stability}, where $\lyacf$ quantifies the level of relaxation. Setting $\lyacf = 0$ recovers the standard CLF \cite{rawlings2017model, grimm2005model}. Moreover, $\lyacf$ is allowed to be $+\infty$, meaning \eqref{eq:sec_prob:relaxed_clf} can always be satisfied. Taking $\ell_x$ as an example, the widely used quadratic cost \cite{eren2017model, rawlings2017model, moreno2022performance, liu2024stability} satisfies Assumption \ref{ass:sec_prob:terminal_cost}. 
However, costs using the $1$-norm or $\infty$-norm, though also used in MPC \cite[Chapter 12.5]{borrelli2017predictive}, do not satisfy Assumption \ref{ass:sec_prob:terminal_cost}.

\begin{remark}[Terminal Cost]
	Choosing $\cost^\star$ as a virtual terminal cost in $\probmpc{\para}$ is equivalent to having no terminal costs \cite{grune2008infinite, grune2017nonlinear, kohler2023stability}, and thus a longer prediction horizon may be required to ensure stability \cite{grune2017nonlinear}. However, designing non-trivial terminal costs usually requires knowing the true model \cite{rawlings2017model, kohler2021stability}, and its analysis under model mismatch is an open problem. For the CLF property \eqref{eq:sec_prob:relaxed_clf}, standard results assumes access to the true model \cite{kohler2021stability, mayne2000constrained}. In this work, we assume that the relaxed CLF condition holds on $\xroa$, following the setting in \cite[Assumption 5]{kohler2023stability}, but only for the nominal model.
\end{remark}

We denote the solution to $\probmpc{\para}$ by $\bu^\star_{N}(x;\para) = \{u^\star_{k|t}(x;\para)\}^{N-1}_{k=0}$ and $J_{N}(x;\bu^\star_{N}(x_t;\para), \para)$ by $V_{N}(x;\para)$. In addition, $\otx{\para}(k;x,\bu^\star_{k}(x;\para))$ is denoted by $\xi^\star_{k}(x;\para)$ for simplicity. Assumptions \ref{ass:sec_pre:basic_problem} and \ref{ass:sec_prob:terminal_cost} ensures $V_{N}(\cdot;\para): \cX \to \bR_+$ is continuous \cite[Theorem~2.7]{rawlings2017model}. Solving $\mathrm{P}_{\mathrm{MPC}}(\theta)$ implicitly defines an MPC policy $\mu_{N,\para}: \cX \to \cU$ as $\mu_{N,\para}(x) = u^\star_{0|t}(x;\para)$, under which the infinite-horizon performance is $J_{\infty}(x;\mu_{N,\para}, \para) :=  \sum^{\infty}_{k = 0}\cost\left(\ctx{\theta}(k;x,\mu_{N,\para}), \ctu{\theta}(k;x,\mu_{N,\para})\right)$. In practice, obtaining the true model is a stringent requirement, and thus we can only apply MPC using an \textit{estimated} nominal model, aligning with the certainty-equivalence principle \cite{hespanha1999certainty}. To guarantee the existence of a well-defined\footnote{For an optimization problem $\min_{z \in \mathcal{Z}}c(z)$ where $c \not\equiv \infty$, existence of a minimizer $z^\star \in \mathcal{Z}$ ensures that $c(z^\star) < \infty$.} CE-MPC policy $\mu_{N,\hpara}$, we impose the following assumption:
\begin{assumption}[Cost Controllability]
	\label{ass:sec_prob:nominal_cost_controllability}
	For all $x \in \xroa$, there exist constants $\umpc < +\infty$ and $\mpc_{i} \in (0, \umpc]$ ($i \in \bI_{1:\infty}$), such that
	\begin{equation}
		\label{eq:sec_prob:mpc_linear_bound}
		V_i(x; \hpara) \leq (1 + \mpc_{i}) \cost^\star(x), \quad \forall i \in \bI_{1:\infty}.
	\end{equation}
\end{assumption}

Assumption \ref{ass:sec_prob:nominal_cost_controllability}, being standard in the literature on RDP \cite[Assumption. 4.2]{grune2008infinite}, \cite[Assumption. 4]{kohler2023stability}, provides a linear bound on $V_i(x; \hpara)$ in terms of $\cost^\star$. This bound connects the infinite-horizon performance to the prediction horizon \cite{grune2017nonlinear}. In the current work, we require \eqref{eq:sec_prob:mpc_linear_bound} to hold for the nominal model instead of the true model as in \cite{grune2017nonlinear}. Besides, Assumption \ref{ass:sec_prob:nominal_cost_controllability} implies that $V_{\infty}(x; \hpara)$ is finite for any fixed $x$, and $\mpc_{i}$ can be approximated numerically by open-loop simulations \cite{kohler2023stability, grune2017nonlinear}. In fact, even though $\mu_{N,\hpara}$ is used, the system evolves based on the true model, so the closed-loop system is

\vspace{-0.25cm}  
\begin{equation}
	\label{eq:sec_prob:closed_loop_system}
	x_{t+1} = f(x_t, \mu_{N,\hpara}(x_t); \tpara).
\end{equation}
Accordingly, the performance of CE-MPC is $J_{\infty}(x;\mu_{N,\hpara}, \tpara)$, and the objectives are specifying under which conditions $ \mu_{N,\hpara}$ stabilizes the true model and analyzing the \textit{suboptimality} of $J_{\infty}(x;\mu_{N,\hpara}, \tpara)$ with respect to $V_{\infty}(x;\tpara)$.

\vspace*{-0.1cm}
\section{Perturbation Analysis of MPC Value Function}
\label{sec:theoretical_analysis}
This section provides the perturbation analysis of the MPC value function, quantifying the variation of $V_N(x; \para)$ when perturbing, respectively, $x$ and $\para$. Notably, we establish a novel property---\textit{consistent error matching} (see Proposition \ref{prop:error_matching})---to relate the variation to $\ell^\ast(x)$. In addition, the asymptotic bounds of core functions in the final results are characterized as side results, considering $\mpara \to 0$ and $N \to \infty$.

In Section \ref{subsec_thereticalAnalysis:preparatory}, we review existing results on sensitivity analysis. Next, variation of $V_N(x; \para)$ due to one-step-ahead prediction error is analyzed in Section \ref{subsec_theoreticalAnalysis:state_perturbation_mpc}. Finally, Section \ref{subsec_theoreticalAnalysis:parameter_perturbation_mpc} focuses on the perturbation with respect to $\para$.
\subsection{Preparatory Analysis}
\label{subsec_thereticalAnalysis:preparatory}
To establish the perturbation bounds, sensitivity analysis for optimal control is needed. Before presenting the main result, we first introduce the following assumption:
\begin{assumption}
	\label{ass:sec_analysis:ssosc_licq}
	$\forall x \in \cX$, $\probmpc{\para}$, together with $\epara$ in \eqref{eq:sec_pre:parameter_set} satisfies the following conditions:
	\begin{enumerate}
		\item (Linear Independence Constraint Qualification (LICQ) (cf. \cite{robinson1980strongly}, \cite[Assumption H.1]{lin2022bounded})) the active-constraint Jacobian\footnote{The active-constraint Jacobian is the Jacobian of all equality constraints and active inequality constraints.}at the optimal primal-dual solution\footnote{The union of the primal variables and the dual variables at the optimum of a constrained optimization problem.} has full row rank for all $\para \in \setpara$;
		\item (Strong Second-Order Sufficient Condition (SSOSC) (cf. \cite{robinson1980strongly}, \cite[Assumption H.1]{lin2022bounded})) the reduced Hessian\footnote{The reduced Hessian is defined as $Z^\top H Z$, where $H$ is the Hessian of the Lagrangian and $Z$ is the null space matrix of the active-constraint Jacobian.}at the optimal primal-dual solution of the problem $\probmpc{\para}$ is positive definite for all $\para \in \setpara$;
		\item (Uniform Singular Spectrum Bounds (cf. \cite[Theorem 4.5]{shin2022exponential} and \cite[Defintion H.2]{lin2022bounded})) There exist positive singular spectrum bounds $\ussbh$, $\ussbr$, and $\lssbh$ for the optimization problem $\probmpc{\para}$ specified for all $\para \in \setpara$, and all $x' \in \mathcal{B}(x, R(x;\epara))$.
	\end{enumerate}
\end{assumption}
The SSOSC and LICQ conditions are standard in sensitivity analysis of general optimization problems \cite[Chapter 4]{bonnans2013perturbation},\cite[Assumption 4.2]{shin2022exponential} and MPC \cite[Chapter 8.6.1]{rawlings2017model}, \cite[Theorem III.3]{wabersich2022cautious}, and they ensure that $\bu^\star_{N}(x_t;\para)$ is a regular minimizer\footnote{For $\min_{z \in \mathcal{Z}}c(z)$, a regular minimizer $z^\star$ is a strict KKT point, i.e., it, together with its Lagrangian duals, satisfies the KKT condition, and there exists $\delta > 0$ such that $f(z') > f(z^\star)$, $\forall z' \in \mathcal{B}(z^\star,\delta) \cap \mathcal{Z}$.}, and thus being locally Lipschitz. To analyze suboptimality, Lipschitzness is frequently assumed in learning-based MPC \cite{liu2025regret, dogan2023regret}, certainty-equivalence MPC \cite{lin2021perturbation, lin2022bounded}, and approximate MPC \cite{alsmeier2024imitation}. Besides, other similar alternatives have been used to ensure regularity properties \cite{mestres2025regularity}, and real-world examples admitting a Lipschitz MPC policy (e.g., inventory control \cite{lin2022bounded}, robot arms \cite{li2025learning}, and drones \cite{wabersich2022cautious}) have also been considered. The third condition means that the singular values of the Hessian of the Lagrangian derived from the MPC optimization are uniformly bounded, which holds when the Lagrangian admits a bounded second derivative and is strongly convex on a bounded domain \cite[Chapter 2]{nesterov2013introductory}. This assumption is useful to establish decaying perturbation bounds in MPC when suboptimality analysis is needed \cite{lin2022bounded, li2025learning}. For MPC, when the constraints are inactive, spectrum bounds can be alternatively checked via controllability and detectability \cite{shin2021controllability}. In general, theoretically identifying the uniform singular spectrum bounds is hard \cite{xu2018exponentially}, and it is not applicable to all systems. However, if the linearized dynamics and costs satisfy the SSOSC and LICQ conditions for the active constraints, then singular spectrum bounds can also be justified \cite{shin2022exponential}.
\begin{lemma}[Exponential Decay of Sensitivity \cite{lin2022bounded, shin2022exponential}] 
\label{lm:sec_analysis:eds}
Let Assumptions \ref{ass:sec_pre:basic_problem}, \ref{ass:sec_prob:terminal_cost}, and \ref{ass:sec_analysis:ssosc_licq} hold. Then there exist constants $\edsc \in (0, \infty)$ and $\edsr \in (0, 1)$ satisfying the following conditions:
\begin{enumerate}
	\item For all $\para \in \setpara$ and any $x \in \xroa$, it holds that
	\vspace*{-0.15cm}
	\begin{equation}
		\label{eq:subsec_eds:parametric_sensitivity}
		\hspace*{-3ex} \|u^\star_{k}(x;\para) \sminus u^\star_{k}(x;\hpara)\| \sleq \edsc\Lambda_N(k) \mpara,\; \forall k \in \bI_{0:N-1},
	\end{equation}
	where $\Lambda_N(k) = \left(\sum^{k}_{i=0}\edsr^i + \sum^{N - k-1}_{i=1}\edsr^i\right)$.
	\item For any $x' \in \cX$ and $x'' \in \mathcal{B}(x', R(x';\epara))$, it holds that
	\vspace*{-0.15cm}
	\begin{equation}
		\label{eq:subsec_eds:state_sensitivity_input}
		\hspace*{-5ex} \|u^\star_{k}(x';\para) \sminus u^\star_{k}(x'';\para)\| \sleq \edsc \edsr^k\|x' \sminus x''\|, \forall k \in \bI_{0:N-1}.
	\end{equation}
\end{enumerate}
The constant $\edsc$ and decay factor $\edsr$ are given, respectively, by
\begin{equation*}
	\edsc = \left(\frac{\ussbh \ussbr}{\lssbh^2}\right)^{\frac{1}{2}},\quad \text{and }
	\edsr = \left(\frac{\ussbh^2 - \lssbh^2}{\ussbh^2 + \lssbh^2}\right)^{\frac{1}{8}}.
\end{equation*} 
\end{lemma}

The constants $H$ and $\gamma$ can be theoretically characterized in special cases \cite{lin2022bounded}, but mostly are determined numerically \cite{shin2022exponential, tranossub, li2025learning}. It is obvious that $\Lambda_N(k) \leq \frac{2}{1 - \edsr}$ for all $k \in \bI_{0:N-1}$ and $N \in \bI_{1:\infty}$. The constant $\Lambda_N(k)$ will be frequently used to simplify the analysis. Apart from sensitivity bounds, additional assumptions on $f$, $\cost_x$, and $\cost_u$ are needed to quantify the perturbation of the MPC value function.
\begin{assumption}[Lipschitz Continuous Dynamics]
	\label{ass:sec_analysis:Lipschitz_dynamics}
	For $x', x'' \in \cX$ and $u', u'' \in \cU$, there exist functions $L_{f,x}, L_{f,u}: \bR^{n_{\para}} \to \bR_{+}$ such that
	\vspace*{-0.25cm}
	\begin{multline}
		\label{eq:lipschitz_dynamics}
		\|f(x', u';\para) - f(x'', u'';\para)\| \leq L_{f,x}(\para)\|x' - x''\| \\ + L_{f,u}(\para)\|u' - u''\|, \;\forall \para \in \setpara,
	\end{multline}
	where $L_{f,x}(\para)$ and $L_{f,u}(\para)$ are the local Lipschitz constants for a given $\para$. In addition, there exist uniform Lipschitz constants $\blfx := \sup_{\para \in \setpara} L_{f,x}(\para) < \infty$ and $\blfu := \sup_{\para \in \setpara} L_{f,u}(\para) < \infty$.
\end{assumption}
Lipschitz continuous systems are widely considered in optimal control \cite{wagenmaker2023optimal, karapetyan2025closed} and MPC \cite{wabersich2022cautious, lin2022bounded, dogan2023regret}. Exact Lipschitz constants are mostly approximated numerically \cite{cobzacs2019lipschitz}. 
\begin{remark}
	For $N < +\infty$ and $\epara < +\infty$, the state and the worst-case state-trajectory perturbation due to mismatch are both bounded, leading to the existence of a \textit{compact} set $\cX_{\text{cp}} \subset \cX$ that includes all possible state trajectories. In this case, Assumption \ref{ass:sec_analysis:Lipschitz_dynamics} becomes redundant since Assumption \ref{ass:sec_pre:basic_problem}-(1) implies \textit{local} Lipschitzness on $\cX_{\text{cp}} \times \cU$.
\end{remark}
\begin{lemma}[Model Difference Bound]
	\label{lm:dynamic_difference}
	Let Assumption \ref{ass:sec_analysis:Lipschitz_dynamics} hold. There exist positive definite functions $\lddx$ and $\lddu$ that are strictly increasing, such that
	\begin{equation}
		\label{eq:difference_lipschitz_dynamics}
		\|\dsys\| \leq \lddx(\mpara)\|x\| + \lddu(\mpara)\|u\|.
	\end{equation}
\end{lemma}
The condition \eqref{eq:difference_lipschitz_dynamics} is a weak extension of \eqref{eq:lipschitz_dynamics}, and the proof of Lemma \ref{lm:dynamic_difference} is given in Appendix \ref{appendix:A0}. To ease the following analysis, we further define the \textit{mismatch Lipschitz function} as
\begin{equation}
	\label{eq:big_O_ease_difference_upper_bound}
	\ldd(\cdot) := \max\{\lddx(\cdot), \lddu(\cdot)\},
\end{equation}
which quantifies the growth behavior of $\|\dsys\|$ with respect to $\mpara$. For the cost functions $\ell_x$ and $\ell_u$, unlike most of the existing results where Lipschitz continuity is employed \cite{schwenkel2020robust, zavala2009advanced,wabersich2022cautious, dogan2023regret, agarwal2019online} such that a \textit{linear} perturbation bound is valid, Assumption \ref{ass:sec_prob:terminal_cost} provides refined bound:
\begin{corollary}[Cost Bounds]
	\label{coro:quadratic}
	Let Assumption \ref{ass:sec_prob:terminal_cost} hold. Then $\ell_x$ and $\ell_u$ satisfy the following conditions:
	\begin{enumerate}
		\item (Perturbation Bound)
		\begin{subequations}
			\label{eq:cost_perturbation}
			\begin{align}
				\label{eq:cost_perturbation_x}
				& \hspace*{-3ex} |\ell_x(x \splus \delta_x) - \ell_x(x)| \leq \frac{L_{\ell,x}}{2}\|\delta_x\|^2 + L_{\ell,x}\|x\|\|\delta_x\|, \\
				\label{eq:cost_perturbation_u}
				& \hspace*{-3ex} |\ell_u(u \splus \delta_u) - \ell_u(u)| \leq \frac{L_{\ell,u}}{2}\|\delta_u\|^2 + L_{\ell,u}\|u\|\|\delta_u\|;
			\end{align}
		\end{subequations}
		\item (Quadratic Lower Bound)
		\begin{equation}
			\label{eq:quadratic_lower}
			\ell_x(x) \geq \frac{m_{\ell,x}}{2}\|x\|^2,\;
			\ell_u(u) \geq \frac{m_{\ell,u}}{2}\|u\|^2.
		\end{equation}
	\end{enumerate}
\end{corollary}

The results \eqref{eq:cost_perturbation} and \eqref{eq:quadratic_lower} can be verified, respectively, using strong smoothness and strong convexity. Corollary \ref{coro:quadratic} helps evaluate the perturbation of $V_N(\cdot;\para)$, and a similar quadratic perturbation bound \qt{will} then follow as well. To adopt the RDP analysis, we establish an instrumental result , called \textit{consistent error matching}, which is given as follows:
\begin{proposition}[Consistent Error Matching]
	\label{prop:error_matching}
	Let Assumptions \ref{ass:sec_prob:terminal_cost} and \ref{ass:sec_analysis:Lipschitz_dynamics} hold. Then the model error \eqref{eq:dynamic_parametric_difference} satisfies
	\begin{equation}
		\label{eq:error_consistency}
		\|\dsys\|^2 \leq \mec \ldd^2(\mpara)\ell(x, u),
	\end{equation}
	where $\mec = 2(m^{-1}_{\ell,x} + m^{-1}_{\ell,u})$ and $\ldd$ is defined in \eqref{eq:big_O_ease_difference_upper_bound}.
\end{proposition}

The proof of Proposition \ref{prop:error_matching} is given in Appendix \ref{appendix:A2}. \eqref{eq:error_consistency} provides an upper bound for the model error in terms of the stage cost $\ell$, which is essential to establish the performance bounds using RDP. Besides, the bound is consistent in that $\|\dsys\| = \ldd(0) = 0$ when $\mpara = 0$. This consistent error bound \eqref{eq:error_consistency} is essential to studying the stability of CE-MPC as well as establishing consistent performance bounds. To convenience, we further define a shorthand $\xp{\para}$ as
\begin{equation}
	\label{eq:state_next}
	\xp{\para} := f(x, \mu_{N,\hpara}(x);\para), \\
\end{equation}
which stands for the next-step state evolving under the system \eqref{eq:sec_pre:nonlinear_system} parameterized by $\para$ under the CE-MPC policy $\mu_{N,\hpara}$.

\vspace*{-0.3cm}
\subsection{One-step-ahead State Perturbation}
\label{subsec_theoreticalAnalysis:state_perturbation_mpc}
We first present a refined result on nominal stability.
\begin{proposition}[Nominal System Stability]
	\label{prop:nominal_energy_decreasing}
	Let Assumptions \ref{ass:pre:initial_state_roa}, \ref{ass:sec_prob:terminal_cost}, and \ref{ass:sec_prob:nominal_cost_controllability} hold. Then the MPC value function $V_{N}(\cdot;\hpara)$ satisfies
	\vspace*{-0.45cm}
	\begin{subequations}
		\begin{align}
			\label{eq:nominal_bound}
			& \cost^\star(x) \leq V_{N}(x;\hpara) \leq (1 + \mpc_{N})\cost^\star(x) \\
			\label{eq:nominal_decreasing}
			& V_{N}(\xp{\hpara};\hpara) \sminus V_{N}(x;\hpara) \leq \sminus(1 \sminus \epsilon_N)\cost(x, \mu_{N,\hpara}(x)),
		\end{align}
	\end{subequations}
	where $\epsilon_N := \frac{(1 + \umpc)\mpc_{N}\lyacf}{(N-1)\lyacf + N + \umpc}$. In addition, if $N \geq \underline{N}_{0} := 1 + \lceil \frac{(1+\umpc)\mpc_{N}\lyacf-\umpc-1}{1 + \lyacf}\rceil$, then $\epsilon_N \in [0, 1)$, meaning that the nominal system is stable under the nominal MPC control policy $\mu_{N,\hpara}$. Moreover, it holds that
	\begin{equation}
		\label{eq:next_state_value_upper_bound}
		V_{N}(\xp{\hpara};\hpara) \leq (\mpc_N + \epsilon_N)l(x, \mu_{N,\hpara}(x)),
	\end{equation}
	where $\mpc_N$ is defined in Assumption \ref{ass:sec_prob:nominal_cost_controllability}.
\end{proposition}

Proposition \ref{prop:nominal_energy_decreasing} states that $V_{N}(\cdot;\hpara)$ is a Lyapunov function for the nominal model under $\mu_{N,\hpara}$, and its proof is given in Appendix \ref{appendix:B1}. Furthermore, \eqref{eq:nominal_decreasing} quantifies the decrease in energy using $\cost(x, \mu_{N,\hpara}(x))$, which will be used to establish state perturbation of $V_N$. It is also noted that 
\begin{equation*}
	\epsilon_N = \bigo{N^{-1}} \text{, when }N \to \infty. 
\end{equation*}

\begin{remark}[Nominal Stability]
	\label{rmk:nominal_stability}
	Results on the nominal stability without terminal conditions are prevalent \cite[Chapter 5]{grune2017nonlinear}, \cite{kohler2023stability, schwenkel2020robust}. The inequality \eqref{eq:nominal_decreasing} is similar to the one in \cite[Theorems 1 and 5]{kohler2023stability}, but our coefficient $\epsilon_N$ is tighter in terms of the required minimal prediction horizon for stability, thus being less restrictive. Moreover, when $\lyacf = +\infty$, $\epsilon_N$ equals $\epsilon^\star_N := \frac{(1+\umpc)\mpc_{N}}{N-1}$ and $\underline{N}_{0}$ degenerates to $\underline{N}^\star_{0} := 1 + \lceil(1+\umpc)\mpc_{N}\rceil$. If a locally exponentially stabilizing controller is available, a tighter $\epsilon_N$ can be obtained such that $\epsilon_N = \bigo{r^N}$ for some $r \in (0,1)$ when $N \to \infty$ \cite{kohler2021stability, liu2024stability}.
\end{remark}
To quantify the difference between the MPC value function $V_{N}(\cdot;\hpara)$ evaluated, respectively, at $\xp{\hpara}$ and $\xp{\para}$, a perturbation bound of the open-loop state is needed.
\begin{lemma}[Open-loop State Perturbation]
	\label{lm:open_loop_state_perturbation}
	Let Assumptions \ref{ass:sec_pre:basic_problem}, \ref{ass:sec_prob:terminal_cost}, \ref{ass:sec_analysis:ssosc_licq}, and \ref{ass:sec_analysis:Lipschitz_dynamics} hold. Then, for all $\para \in \setpara$, any $N \in \bI_{1:\infty}$, $x' \in \cX$, and $x'' \in \mathcal{B}(x', R_{x'})$, the open-loop state $\xi^\star_{k,N}(\cdot;\para)$ satisfies the following Lipschitz-type perturbation bound:
	\begin{equation}
		\label{eq:lm4_open_loop_state_perturbation_bound}
		\|\xi^\star_{k}(x';\para) - \xi^\star_{k}(x'';\para)\| \leq \Gamma_{k}\|x' - x''\|, \forall k \in \bI_{0:N},
	\end{equation}
	where $\Gamma_{k} := (\blfx)^k + \blfu\edsc\sum^{k-1}_{i=0}(\blfx)^i\edsr^{k-i}$ with $\edsc$ and $\edsr$ given in Lemma \ref{lm:sec_analysis:eds}.
\end{lemma}

\begin{proposition}[One-step-ahead State Perturbation]
	\label{prop:parametric_state_perturbation}
	Let Assumptions \ref{ass:sec_pre:basic_problem} to \ref{ass:sec_analysis:Lipschitz_dynamics} hold. There exists a function $\alpha^\ast_N \in \mathcal{K}$ such that, for all $\para \in \setpara$ and all $x \in \cX$, it holds that
	\begin{equation}
		\label{eq:parametric_state_perturbation}
		\hspace*{-2.5ex} |V_{N}(\xp{\para};\hpara) \sminus V_{N}(\xp{\hpara};\hpara)| \sleq \alpha^\ast_N(\mpara)\cost(x, \mu_{N,\hpara}(x)).\hspace*{-2ex}
	\end{equation}
	In addition, $\alpha^\ast_N(\cdot)$ admits a quadratic form as
	\begin{equation}
		\label{eq:alpha_detail_in_the_theorem}
		\vspace*{-1ex}
		\alpha^\ast_N(\cdot) = \pi_{\alpha,N,2}\ldd^2(\cdot) + \pi_{\alpha,N,1}\ldd(\cdot),
	\end{equation}
	where the details of $\pi_{\alpha,N,2}$ and $\pi_{\alpha,N,1}$ are given, respectively, in \eqref{eq:B3_second_order} and \eqref{eq:B3_first_order} in Appendix \ref{appendix:B3}. Moreover, the asymptotic bounds of $\alpha^\ast_N(\mpara)$ are given as follows:
	\begin{itemize}
		\item $\bigo{\ldd(\mpara)}$ for $\blfx < 1$;
		\item $\bigo{N\ldd(\mpara)}$ for $\blfx = 1$;
		\item $\bigo{\blfx^{2N}\ldd(\mpara)}$ for $\blfx > 1$,
	\end{itemize}
	when $\mpara \to 0$ and $N \to \infty$.
\end{proposition}

Eq. \eqref{eq:lm4_open_loop_state_perturbation_bound} provides a Lipschitz bound of the open-loop state, where $\Gamma_k$ is the Lipschitz constant. Based on Lemma \ref{lm:open_loop_state_perturbation}, we can establish the desired bound in Proposition \ref{prop:parametric_state_perturbation} using $\mpara$. It should be noted that \eqref{eq:parametric_state_perturbation} is \textit{consistent} since $\alpha_N$ is positive definite, i.e., given that $\mpara = 0$, $\xp{\para} = \xp{\hpara}$ and $\alpha^\ast_N(0) = 0$. In general, $\alpha^\ast_N$ inherits the property of the mismatch function \eqref{eq:big_O_ease_difference_upper_bound}. Moreover, if $\ldd$ is linear, then $\alpha^\ast_N(\mpara)$ scales linearly with $\mpara$ for an infinitesimal mismatch.
\begin{remark}[Open-loop System \& Contraction]
	\label{rmk:open_loop_contraction}
	Lemma \ref{lm:open_loop_state_perturbation} is coarse since only the Lipschitz continuity is used without other properties of the open-loop system, and thus $\Gamma_k$ exhibits exponential growth as $k$ increases, making \eqref{eq:parametric_state_perturbation} also \textit{conservative}, especially for $\overline{L}_{f,x} \geq 1$. Nonetheless, refinements are possible when adding terminal constraints \cite{lin2022bounded}, assuming the existence of an exponentially stable terminal policy \cite{kohler2021stability}, or considering pre-stabilization such that $\overline{L}_{f,x} < 1$ or the system satisfies other contraction properties \cite{tran2018convergence, karapetyan2025closed}. However, under input constraints, achieving such contraction behaviors is not always possible, especially when the system exhibits rapidly growing dynamics (e.g., high-order polynomials and exponential functions). These potential refinements to achieve contraction and to obtain tighter, contractive bounds are not essential for this paper and are left for future work.
\end{remark}

\vspace*{-0.5cm}
\subsection{Parameter Perturbation}
\label{subsec_theoreticalAnalysis:parameter_perturbation_mpc}
Next, we investigate the difference between the values $V_{N}(x;\hpara)$ and $V_{N}(x;\tpara)$. Noting that $V_{N}(\mathbf{0};\hpara) = V_{N}(\mathbf{0};\tpara) = 0$ under Assumptions \ref{ass:sec_pre:basic_problem} and \ref{ass:sec_prob:terminal_cost}, it is desired to obtain a bound on this difference in terms of $\cost^\star$ such that it scales with $x$. We first give a scalable input sensitivity bound as follows:

\begin{lemma}[Scalable Input Perturbation Bound]
	\label{lm:scalable_eds}
	Let Assumptions \ref{ass:sec_pre:basic_problem}, \ref{ass:sec_prob:terminal_cost}, and \ref{ass:sec_analysis:ssosc_licq} hold. Define functions $\sueds{k}(\cdot;d_x): \bR_+ \to \bR_+$ for all $k \in \bI_{0:N-1}$ by
	\begin{multline*}
		\sueds{k}(\delta_{\para};d_x) := \\
		\begin{cases}
			\edsc\Lambda_N(k)R(\mathbf{0};\epara)^{-1}\delta_{\para}d_x & x \in \cX \setminus \mathcal{B}(\mathbf{0},R(\mathbf{0};\epara)) \\
			\min\{\Lambda_N(k)\delta_{\para}, 2\edsr^k d_x\}\edsc & x \in \mathcal{B}(\mathbf{0},R(\mathbf{0};\epara))
		\end{cases},
	\end{multline*}
	where $\edsc$ and $\Lambda_N(k)$ are defined as in Lemma \ref{lm:sec_analysis:eds}. Then, for all $\para \in \setpara$ and $x \in \cX$, the following holds:
	\begin{equation}
		\label{eq:subsec_para_perturbation:refined_eds}
		\hspace*{-1ex}\|u^\star_{k}(x;\hpara) - u^\star_{k}(x;\para)\| \leq \sueds{k}(\mpara; \|x\|),\;\forall k \in \bI_{0:N-1}.
	\end{equation}
	By definition, it holds that, for all $k \in \bI_{0:N-1}$, $\sueds{k}(0;d_x) = 0$, and $\sueds{k}(\cdot, d_x)$ is non-decreasing. Besides, there exists a function $\overline{\omega}(\cdot;d_x): \bR_+ \to \bR_+$ defined by
	\begin{equation*}
		\overline{\omega}(\delta_{\para};d_x) := 
		\begin{cases}
			\frac{2\edsc}{(1-\edsr)R(\mathbf{0};\epara)}\delta_{\para}d_x & x \in \cX \setminus \mathcal{B}(\mathbf{0},R(\mathbf{0};\epara)) \\
			2\min\{\frac{\delta_{\para}}{1-\edsr}, d_x\}\edsc & x \in \mathcal{B}(\mathbf{0},R(\mathbf{0};\epara))
		\end{cases}
	\end{equation*}
	such that $\sueds{k} \leq \overline{\omega}$ for all $k \in \bI_{0:N-1}$ and $N \in \bI_{1:\infty}$.
\end{lemma}

\begin{lemma}[Open-loop Parametric Perturbation]
	\label{lm:open_loop_parameter_perturbation}
	Let Assumptions \ref{ass:sec_pre:basic_problem}, \ref{ass:sec_prob:terminal_cost}, \ref{ass:sec_analysis:ssosc_licq}, and \ref{ass:sec_analysis:Lipschitz_dynamics} hold. Then, for all $\para \in \setpara$, $N \in \bI_{1:\infty}$, and any $x \in \cX$, the open-loop state $\xi^\star_{k}(\cdot;\para)$ satisfies the following accumulative perturbation bound for all $k \in \bI_{0:N}$:
	\begin{multline}
		\label{eq:open_loop_parameter_perturbation_bound}
		\|\xi^\star_{k}(x;\para) - \xi^\star_{k}(x;\hpara)\| \leq C_{\Delta u,k}(\mpara;x) + \\ \sum^{k-1}_{i=0}\blfx^{k-1-i}\|\Delta f(\xi^\star_{i}(x;\hpara),u^\star_{i}(x;\hpara);\para)\|,
	\end{multline}
	where $C_{\Delta u,k}(\cdot;x) := \blfu$$\sum^{k-1}_{i=0}\blfx^{k-1-i}\sueds{i}(\cdot;\|x\|)$ with $\sueds{i}$ defined as in Lemma \ref{lm:scalable_eds}. 
\end{lemma}
Lemma \ref{lm:scalable_eds} helps establish the parametric perturbation of the open-loop state in Lemma \ref{lm:open_loop_parameter_perturbation}, and the corresponding parameter perturbation of $V_N(x, \hpara)$ is summarized in Proposition \ref{prop:general_parameter_perturbation}.

\begin{proposition}[General Affine Parameter Perturbation]
	\label{prop:general_parameter_perturbation}
	Let Assumptions \ref{ass:sec_pre:basic_problem} to \ref{ass:sec_analysis:Lipschitz_dynamics} hold. Then for all $N \in \bI_{1:\infty}$, $\para \in \setpara$, and $x \in \cX$, it holds that
	\begin{equation}
		\label{eq:parameter_perturbation}
		|V_{N}(x;\hpara) - V_{N}(x;\para)| \leq \beta_N(\mpara; x),
	\end{equation}
	where $\beta_N(\cdot; x)$ is non-decreasing satisfying $\beta_N(0; x) = 0$, and it admits a simple form as 
	\begin{multline*}
		\beta_N(\cdot; x) = \pi_{\beta,N,2}(x)\ldd^2(\cdot) + \pi_{\beta,N,1}(x)\ldd(\cdot) + \\ \zeta_{\beta,N,2}[\overline{\omega}(\cdot;\|x\|)]^2 + \zeta_{\beta,N,1}(x)\overline{\omega}(\cdot;\|x\|),
	\end{multline*}
	and the details of $\pi_{\beta,N,2}(x)$, $\pi_{\beta,N,1}(x)$, $\zeta_{\beta,N,2}$, and $\zeta_{\beta,N,1}(x)$ are given in \eqref{eq:C3_final} in Appendix \ref{appendix:C3}. Moreover, the asymptotic bounds of $\beta_N(\mpara; x)$ are given as follows:
	\begin{itemize}
		\item $\bigo{N(\ldd(\mpara) + \mpara)}$ for $\blfx < 1$;
		\item $\bigo{N^3(\ldd(\mpara) + \mpara)}$ for $\blfx = 1$;
		\item $\bigo{\blfx^{2N}(\ldd(\mpara) + \mpara)}$ for $\blfx > 1$,
	\end{itemize}
	when $\mpara \to 0$ and $N \to \infty$.
\end{proposition}
The perturbation function $\beta_N(\cdot; x)$ in Proposition \ref{prop:general_parameter_perturbation} depends on both $\ldd(\cdot)$ and $\overline{\omega}(\cdot;\|x\|)$ in general. However, under a mild assumption of the input sensitivity, a refined bound can be obtained such that it only depends on $V_{N}(x;\para)$.
\begin{assumption}[Local Separable Input Perturbation]
	\label{ass:local_scalable_input_perturbation_bound}
	There exist a neighborhood $\Omega \subseteq \mathcal{B}(\mathbf{0}, R(\mathbf{0};\epara))$ of the origin and constants $\{\eta_{k,N}\}^{N-1}_{k=0}$, such that
	\begin{equation}
		\label{eq:local_scalable_input_perturbation_bound}
		\|u^\star_{k}(x;\hpara) - u^\star_{k}(x;\para)\| \leq \eta_{k,N}\mpara\|x\|, \; \forall k \in \bI_{0:N-1},
	\end{equation}
	and $\eta_{k,N} \leq \bar{\eta}$ for all $k \in \bI_{0:N-1}$ and $N \in \bI_{1:\infty}$.
\end{assumption}
\begin{remark}[Separable Input Perturbation]
	\label{rmk:separable_local_input_perturbation}
	Assumption \ref{ass:local_scalable_input_perturbation_bound} offers a perturbation bound that scales with $\|x\|$, which is compatible with $\|u^\star_{k|t}(x;\para)\| \leq \edsc\edsr^k\|x\|$ since $u^\star_{k|t}(\mathbf{0};\para) = 0$ (cf. \eqref{eq:subsec_eds:state_sensitivity_input} in Lemma \ref{lm:sec_analysis:eds}). Besides, given $g_u(\bzero) < \mathbf{1}$, $\|u^\star_{k|t}(x;\para)\| \leq \edsc\edsr^k\|x\|$ implies that a small enough neighborhood $\Omega$ of the origin exists such that the optimization problem $\probmpc{\para}$ is \textit{unconstrained} within $\Omega$, i.e., all input constraints are \textit{inactive}. Thus, $u^\star_{k|t}(\cdot;\para)$ is locally continuous under Assumptions \ref{ass:sec_pre:basic_problem}, \ref{ass:sec_prob:terminal_cost} and \ref{ass:sec_analysis:ssosc_licq}. In this case, if $u^\star_{k|t}(x;\para) = K_{k,N}(\para)\mu_{\Omega}(x)$ (i.e., the control law is separable) and $K_{k,N}$ is Lipschitz continuous (i.e., $\|K_{k,N}(\para) - K_{k,N}(\hpara)\| \leq L_{K,k,N}\mpara$), then establishing $\mu_{\Omega}(x) \leq M_{\Omega}\|x\|$ leads to\footnote{Based on the continuity of $u^\star_{k|t}(\cdot;\para)$, leveraging the Taylor expansion, $\mu_{\Omega}(x) \leq M_{\Omega}\|x\|$ means the dominant term of $\mu_{\Omega}(x)$ within $\Omega$ is linear.} $\eta_{k,N} = M_{\Omega}L_{K,k,N}$ in \eqref{eq:local_scalable_input_perturbation_bound}. The above reasoning highlights that Assumption \ref{ass:local_scalable_input_perturbation_bound} is not restrictive, and it aligns with the EDS property and the considered MPC problem, and such separation structures have also been used in \cite{kuntz2024beyond} for Lyapunov functions. In fact, for linear systems with quadratic costs, Assumption \ref{ass:local_scalable_input_perturbation_bound} holds since the optimal input to the unconstrained linear quadratic control problem is separable and linear in $x$. Formal characterization and computation of $\eta_{k,N}$ for general nonlinear systems is out of the scope of this paper and it requires further investigation for future work.
\end{remark}
Assumption \ref{ass:local_scalable_input_perturbation_bound} leads to Corollaries \ref{coro:scalable_eds} and \ref{coro:open_loop_parameter_perturbation}, which are the refined version of Lemmas \ref{lm:scalable_eds} and \ref{lm:open_loop_parameter_perturbation}, respectively. In addition, a refined perturbation follows directly as in Proposition \ref{prop:linear_parameter_perturbation}, making the final bound \textit{linear} in $V_N(x;\para)$. Besides, if $\ldd$ in \eqref{eq:big_O_ease_difference_upper_bound} is a linear function of $\mpara$, then $\beta^\ast_N(\mpara)$ scales linearly with $\mpara$ for an infinitesimal mismatch.

\begin{corollary}[Separable Scalable Input Perturbation Bound]
	\label{coro:scalable_eds}
	Let Assumptions \ref{ass:sec_pre:basic_problem}, \ref{ass:sec_prob:terminal_cost}, \ref{ass:sec_analysis:ssosc_licq}, and \ref{ass:local_scalable_input_perturbation_bound} hold. Given $N \in \bI_{1:\infty}$, define constants $\{\leds{k}^\ast\}^{N-1}_{k=0}$ as
	\begin{equation*}
		\leds{k}^\ast := \max\left\{\frac{\edsc\Lambda_N(k)}{\inf_{x\in\xroa\setminus\Omega}\|x\|}, \eta_{k, N}\right\},
	\end{equation*}
	where $\Lambda_N(k)$ and $\Omega$ are defined in Lemma \ref{lm:sec_analysis:eds} and Assumption \ref{ass:local_scalable_input_perturbation_bound}, respectively. Then, $\forall \para \in \setpara$ and $x \in \cX$, it holds that
	\begin{equation}
		\label{eq:subsec_para_perturbation:refined_eds_separable}
		\|u^\star_{k}(x;\hpara) - u^\star_{k}(x;\para)\| \leq \leds{k}^\ast \mpara \|x\|, \; \forall k \in \bI_{0:N-1}.
	\end{equation}
	In addition, $\leds{k}^\ast \leq \bar{\eta}^\ast$ for all $k \in \bI_{0:N-1}$ and $N \in \bI_{1:\infty}$, where $\bar{\eta}^\ast = \max\left\{\frac{2\edsc}{(1 - \edsr )\inf_{x\in\xroa\setminus\Omega}\|x\| }, \bar{\eta}\right\}$.
\end{corollary}

\begin{corollary}[Separable Open-loop Parametric Perturbation]
	\label{coro:open_loop_parameter_perturbation}
	Let Assumptions \ref{ass:sec_pre:basic_problem}, \ref{ass:sec_prob:terminal_cost}, \ref{ass:sec_analysis:ssosc_licq}, \ref{ass:sec_analysis:Lipschitz_dynamics}, and and \ref{ass:local_scalable_input_perturbation_bound} hold. Then, for all $\para \in \setpara$, any $N \in \bI_{1:\infty}$, and any $x \in \cX$, the open-loop state $\xi^\star_{k}(\cdot;\para)$ satisfies the following accumulative perturbation bound for all $k \in \bI_{0:N}$
	\begin{multline}
		\label{eq:co4_open_loop_parameter_perturbation_bound}
		\|\xi^\star_{k}(x;\hpara) - \xi^\star_{k}(x;\para)\| \leq C^\ast_{\Delta u,k}\mpara\|x\| + \\ \sum^{k-1}_{i=0}\blfx^{k-1-i}\|\Delta f(\xi^\star_{i}(x;\para),u^\star_{i}(x;\para);\para)\|,
	\end{multline}
	where $C^\ast_{\Delta u,k}:= \blfu$$\sum^{k-1}_{i=0}\blfx^{k-1-i}\leds{i}^\ast$ with $\leds{i}^\ast$ defined in Corollary \ref{coro:scalable_eds}. 
\end{corollary}

\begin{proposition}[Linear Parameter Perturbation]
	\label{prop:linear_parameter_perturbation}
	Let Assumptions \ref{ass:sec_pre:basic_problem} to \ref{ass:local_scalable_input_perturbation_bound} hold. Then for all $N \in \bI_{1:\infty}$, $\para \in \setpara$, and $x \in \cX$, it holds that
	\begin{equation}
		\label{eq:parameter_perturbation_linear}
		|V_{N}(x;\hpara) - V_{N}(x;\para)| \leq \beta^{\ast}_{N}(\mpara)V_N(x;\para),
	\end{equation}
	where $\beta^\ast_N \in \mathcal{K}$, and it admits a simple form as
	\begin{equation*}
		\beta^\ast_N(\cdot) = \pi^\ast_{\beta,N,2}\ldd^2(\cdot) + \pi^\ast_{\beta,N,1}\ldd(\cdot) + \zeta^\ast_{\beta,N,2}(\cdot)^2 + \zeta^\ast_{\beta,N,1}(\cdot),
	\end{equation*}
	where the details of $\pi^\ast_{\beta,N,1}$, $\pi^\ast_{\beta,N,2}$, $\zeta^\ast_{\beta,N,1}$, and $\zeta^\ast_{\beta,N,2}$ are given in \eqref{eq:C4_final} in Appendix \ref{appendix:C4}. Moreover, the asymptotic bounds of $\beta^\ast_N(\mpara)$ are given as follows:
	\begin{itemize}
		\item $\bigo{N(\ldd(\mpara) + \mpara)}$ for $\blfx < 1$;
		\item $\bigo{N^3(\ldd(\mpara) + \mpara)}$ for $\blfx = 1$;
		\item $\bigo{\blfx^{2N}(\ldd(\mpara) + \mpara)}$ for $\blfx > 1$,
	\end{itemize}
	when $\mpara \to 0$ and $N \to \infty$.
\end{proposition}

\begin{remark}[State Constraints, Compounding Error \& Feasibility]
	Compounding error of the state trajectory is a common problem when the input is perturbed \cite{liu2025regret, karapetyan2025closed}. In our case, the input perturbation stems from model mismatch. This compounding error can grow unboundedly, posing a great challenge in theoretical analysis.  When state constraints are present, robust designs (e.g., tube-based MPC \cite{fleming2014robust} and/or constraint tightening \cite{wabersich2022cautious}) are, in general, inevitable, which do not fit the basic CE-MPC framework. For the above reasons, explicit state constraints are not considered in the current work, which also resolves the feasibility issue of MPC.
\end{remark}

\section{Stability and Performance Analysis}
\label{sec:theoreticalAnalysis:stability_suboptimality}
This section presents the joint effect of the parametric model mismatch and the prediction horizon on the stability and performance of the closed-loop system in \eqref{eq:sec_prob:closed_loop_system}. The results build on the perturbation analysis in Propositions \ref{prop:parametric_state_perturbation}, \ref{prop:general_parameter_perturbation}, and \ref{prop:linear_parameter_perturbation}. 


Specifically, a sufficient condition for the stability of CE-MPC is first established, based on which performance bounds of CE-MPC can be derived. Following the stability condition, the maximum allowed mismatch level that guarantees the closed-loop stability can be derived. Moreover, under mild conditions such that a competitive-ratio bound is valid, an approximate optimal prediction horizon that gives the lowest ratio bound can also be numerically computed by minimizing this bound. Finally, we provide a novel analysis pipeline, \qt{sketching} the core steps to establish the results.   

\begin{theorem}[Stability of CE-MPC]
	\label{thm:stability_ce_mpc}
	Let Assumptions \ref{ass:sec_pre:basic_problem} to \ref{ass:sec_analysis:Lipschitz_dynamics} hold. If the prediction horizon $N \geq \underline{N}_{0}$ and the mismatch level $\epara$ in \eqref{eq:sec_pre:parameter_set} satisfy
	\begin{equation}
		\label{eq:stability_ce_mpc}
		\epsilon_N + \alpha_N^\ast(\epara) < 1,
	\end{equation}
	where $\epsilon_N$ is defined in \eqref{eq:nominal_decreasing} and $\alpha_N^\ast(\cdot)$ is given in Proposition \ref{prop:parametric_state_perturbation}, then the closed-loop system \eqref{eq:sec_prob:closed_loop_system} is asymptotically stable at the origin for any $x_0 \in \xroa$ and any possible true system given in \eqref{eq:sec_pre:nonlinear_system} characterized by $\tpara \in \setpara$.
\end{theorem}

Theorem \ref{thm:stability_ce_mpc} provides a sufficient condition for the stability of CE-MPC, and the relationship in \eqref{eq:stability_ce_mpc} reflects that the mismatch level and the prediction horizon jointly influence the stability of the closed-loop system. An alternate approach, in the context of inherent robustness, derives stability of CE-MPC under the assumption of knowing a prior bound of the mismatch level such that a certain decreased energy property holds in the sense of Lyapunov stability, see \cite[Assumption 9]{kuntz2024beyond} and \cite[Theorem 8]{kuntz2024beyond} for more details. In addition, for a given prediction horizon $N \geq \underline{N}_{0}$ (see Proposition \ref{prop:nominal_energy_decreasing}), an upper bound of $\epara$ can be derived from \eqref{eq:stability_ce_mpc}, i.e., 
\begin{equation}
	\label{eq:supremum_mismatch}
	\hspace*{-1.2ex}\epara < L_d^{-1}\left(\frac{-\pi_{\alpha,N,1} + [\pi_{\alpha,N,1}^2 + 4\pi_{\alpha,N,2}(1 - \epsilon_N)]^{\frac{1}{2}}}{2\pi_{\alpha,N,2}} \right),\hspace*{-1ex}
\end{equation}
where $L_d^{-1}(\cdot)$ is the inverse function of $L_d$ defined in \eqref{eq:big_O_ease_difference_upper_bound} and the constants $\pi_{\alpha,1,N}$ and $\pi_{\alpha,2,N}$ are defined in Proposition \ref{prop:parametric_state_perturbation}. Eq. \eqref{eq:supremum_mismatch} provides the maximum allowed mismatch level that guarantees the closed-loop stability of CE-MPC. Based on the insights from Proposition \ref{prop:parametric_state_perturbation}, extending the horizon does not cause an unbounded perturbation when $\blfx < 1$ (i.e., the dynamics itself is state-wise contractive); however, one should carefully restrict the prediction horizon if $\blfx \geq 1$ since $\alpha^\ast_N(\cdot)$ has a polynomial or even exponential growth when $N$ increases, meaning the stability is more vulnerable. 

\begin{remark}[Stability Verification]
	Verifying \eqref{eq:stability_ce_mpc} (or equivalently \eqref{eq:supremum_mismatch}) only requires the information of the nominal model without the knowledge of $\epara$. To elaborate, the upper bound of $\epara$ in \eqref{eq:supremum_mismatch} comprises $\pi_{\alpha,N,1}$ and $\pi_{\alpha,N,2}$, both of which depend on $\Gamma_k$ defined in \eqref{eq:lm4_open_loop_state_perturbation_bound}, thereby being influenced by $\blfx$ and $\blfu$. However, since \eqref{eq:parametric_state_perturbation} is derived for the nominal model, $\Gamma_k$ can instead be defined using the nominal Lipschitz constants $L_{f,x}(\hpara)$ and $L_{f,u}(\hpara)$. Consequently, computing the upper bound does not require the uniform Lipschitz constants, which depend on knowledge of $\epara$ \textit{a priori}. As a result, both \eqref{eq:stability_ce_mpc} and \eqref{eq:supremum_mismatch} remain valid conditions to verify stability. Still, in Lemma \ref{lm:open_loop_state_perturbation} and Proposition \ref{prop:parametric_state_perturbation}, the more conservative uniform Lipschitz constants are used for brevity.
\end{remark}
\begin{remark}[Relation to Inherent Robustness]
	The stability of MPC under additive disturbances and state estimation error has been studied using the notion of inherent robustness \cite{grimm2007nominally, allan2017inherent, pannocchia2011conditions}, revealing that nominal MPC can stabilize the system under sufficiently small perturbations. In contrast, Theorem \ref{thm:stability_ce_mpc}, though conveying the same core message (i.e., stability of the nominal controller is preserved if the mismatch is sufficiently small), is tailored for multiplicative parametric uncertainty without state estimation error, thus being different from the settings in \cite{grimm2007nominally, allan2017inherent, pannocchia2011conditions}. As such, our result can also be viewed as the inherent robustness of MPC under multiplicative uncertainty with a specific stability certificate in terms of mismatch level and the prediction horizon.
\end{remark}

Following Theorem \ref{thm:stability_ce_mpc}, a performance bound characterizing the suboptimality of the CE-MPC controller with respect to the oracle controller derived from $\probopt{\tpara}$ can be obtained, which is given in Theorem \ref{thm:affine_performance_ce_mpc}.

\begin{theorem}[Performance Bounds of CE-MPC]
	\label{thm:affine_performance_ce_mpc}
	Let Assumptions \ref{ass:sec_pre:basic_problem} to \ref{ass:sec_analysis:Lipschitz_dynamics} hold. Besides, assume the prediction horizon $N \geq \underline{N}_{0}$ and the mismatch level $\epara$ in \eqref{eq:sec_pre:parameter_set} satisfy \eqref{eq:stability_ce_mpc}. Then the infinite-horizon performance $J_{\infty}(x;\mu_{N,\hpara}, \tpara)$ satisfies
	\begin{equation}
		\label{eq:performance_ce_mpc_affine}
		\hspace*{-1ex} J_{\infty}(x;\mu_{N,\hpara}, \tpara) \leq \frac{V_{\infty}\big(x;\tpara\big) + \mathcal{G}_N\left(x;\epara\right)}{1 - \big(\epsilon_N + \alpha_N^\ast(\epara)\big)},
	\end{equation}
	where
	\begin{equation}
		\label{eq:performance_ce_mpc_affine_numerator}
		\mathcal{G}_N\left(x;\epara\right) = \min\bigg\{\umpc \cost^\star(x), \beta_N(\epara;x) \bigg\},
	\end{equation}
	in which $\umpc$ is defined in Assumption \ref{ass:sec_prob:nominal_cost_controllability}, $\alpha^\ast_N$ is given in Proposition \ref{prop:parametric_state_perturbation}, and $\beta_N(\epara;x)$ is given in Proposition \ref{prop:general_parameter_perturbation}. In addition, if Assumption \ref{ass:local_scalable_input_perturbation_bound} holds, then
	\begin{equation}
		\label{eq:performance_ce_mpc_linear}
		J_{\infty}(x;\mu_{N,\hpara}, \tpara) \leq \compratio{N}(\epara)V_{\infty}(x;\tpara),
	\end{equation}
	where
	\begin{equation}
		\label{eq:bound_competitive_ratio}
		\compratio{N}(\epara) = \frac{1 + \min\{\umpc, \beta^{\ast}_{N}(\epara)\}}{1 - \big(\epsilon_N + \alpha_N^\ast(\epara)\big)},
	\end{equation}
	in which the function $\beta^{\ast}_{N}(\cdot)$ is defined in Proposition \ref{prop:linear_parameter_perturbation}.
\end{theorem}

\begin{remark}[Conservatism \& Qualitative Value]
	The stability condition \eqref{eq:stability_ce_mpc} is only \textit{sufficient} due to several (inequality) relaxations made to establish $\alpha^\ast_N$ in Proposition \ref{prop:parametric_state_perturbation}. As such, CE-MPC may still stabilize the true system even with a larger mismatch level that invalidates \eqref{eq:stability_ce_mpc}. In addition, since relaxations are also used to derive $\beta_{N}(\cdot;x)$ and $\beta^\ast_N$, the bounds \eqref{eq:performance_ce_mpc_affine} and \eqref{eq:performance_ce_mpc_linear} are both conservative. It is noted that conservatism is a common issue in theoretical performance bounds \cite{grune2017nonlinear, kohler2023stability, lin2022bounded, shi2025suboptimality, giselsson2013feasibility}, especially for systems with \textit{multiplicative} uncertainty where the worst-case prediction error accumulates \textit{exponentially} in general. Therefore, the bounds have limited quantitative value in terms of characterizing the true performance, but are more useful in revealing qualitative trends of the performance with respect to the mismatch level and the horizon. \qt{In general, obtaining sufficiently tight and widely applicable bounds is challenging.}
\end{remark}

Theorem \ref{thm:affine_performance_ce_mpc} provides a general performance bound \eqref{eq:performance_ce_mpc_affine} and a \textit{competitive-ratio} bound \eqref{eq:performance_ce_mpc_linear}. Both of the bounds are expressed using $V_{\infty}(x;\tpara)$, quantifying the suboptimality of the CE-MPC controller compared to the infinite-horizon optimal controller that knows the true model. It should be noted that satisfying the stability condition \eqref{eq:stability_ce_mpc} is a prerequisite to ensure that the bounds are finite and non-negative. Following Propositions \ref{prop:parametric_state_perturbation} and \ref{prop:linear_parameter_perturbation}, asymptotic bounds of $\compratio{N}(\epara)$ in \eqref{eq:bound_competitive_ratio} can also be specified under mild conditions: 
\begin{corollary}
	\label{coro:asymptotic_competitive_ratio}
	Let Assumptions \ref{ass:sec_pre:basic_problem} to \ref{ass:local_scalable_input_perturbation_bound} hold. If $\alpha^\ast_N(\epara) \to 0$ and $\beta^\ast_N(\epara) \to 0$ as $\epara \to 0$ and $N \to \infty$, then the competitive ratio $\compratio{N}(\epara)$ admits the following asymptotic bounds:
	\begin{itemize}
		\item $1 + \mathcal{O}\bigg(\frac{1}{N} + N[\ldd(\epara)+\ldd^2(\epara)] + N\epara[1+\ldd(\epara)]\bigg)$ for $L_{f,x} < 1$;
		\item $1 + \mathcal{O}\bigg(\frac{1}{N} + N^3\ldd(\epara)+N^3\ldd^2(\epara) + N^4\epara + N^4\epara \ldd(\epara)\bigg)$ for $L_{f,x} = 1$;
		\item $1 + \mathcal{O}\bigg(\frac{1}{N} + L^{2N}_{f,x}\ldd(\epara)+L^{2N}_{f,x}\ldd^2(\epara) + L^{4N}_{f,x}\epara+L^{4N}_{f,x}\epara \ldd(\epara)\bigg)$ for $L_{f,x} > 1$,
	\end{itemize}
	when $\epara \to 0$ and $N \to \infty$.
\end{corollary}

For all the cases in Corollary \ref{coro:asymptotic_competitive_ratio}, the only term that makes $\compratio{N}(\epara)$ decrease as $N \to \infty$ is $\frac{1}{N}$ coming from $\epsilon_N$, and this factor can be improved if a local exponentially stabilizing control law exists \cite{kohler2021stability, liu2024stability} (see Remark \ref{rmk:nominal_stability}). The other terms that couple with $\epara$ manifest that increasing $N$ may lead to a larger competitive ratio when model mismatch is present. Such a tradeoff between increasing the horizon to access more future information and restricting it to avoid the accumulation of prediction errors is quantified by the performance bounds \eqref{eq:performance_ce_mpc_affine} and \eqref{eq:performance_ce_mpc_linear}. In fact, given the competitive ratio in \eqref{eq:bound_competitive_ratio}, an \textit{approximate} optimal prediction horizon can be obtained as
\begin{equation}
	\label{eq:optimal_prediciton_horizon}
	N^\star(\epara) = \min_{N \in \bI_{\underline{N}_0:\infty}} \compratio{N}(\epara),
\end{equation}
which can be solved numerically or by enumeration.

\begin{remark}
	The approximate optimal horizon $N^\ast(\epara)$ is a side product of the performance analysis. Given that $N^\star(\epara)$ is derived from the conservative performance bound instead of the true performance, it is, in general, not the real optimal horizon when applying CE-MPC to the true model. Besides, since the performance bounds cannot quantitatively reflect the true performance, using $N^\ast(\epara)$ (or the true optimal horizon) may not exhibit significant performance advantages. Nevertheless, $N^\ast(\epara)$ can be used as a coarse warm start for horizon tuning when lacking perfect system knowledge and/or when conducting real experiments is expensive.
\end{remark}

The analysis pipeline summarizing the different steps and all the important results is provided in Fig. \ref{fig:pipeline}. Moreover, to obtain less conservative performance bounds in cases where extra assumptions are available or the system is subject to \textit{additive} uncertainty, one could refine the nominal stability analysis to obtain a tighter coefficient $\epsilon_N$, as well as polish the perturbation analysis of the MPC value function to establish tighter perturbation functions (i.e., $\alpha^\ast_N$ and $\beta^\ast_{N}$) for the cost functions considered in Assumption \ref{ass:sec_prob:terminal_cost} or other classes of cost functions. With the refined coefficient and/or perturbation functions, other variations of the stability condition and performance bounds follow directly using our analysis pipeline. With the improved performance bound with reduced conservatism, the mismatch upper bound in \eqref{eq:supremum_mismatch} will be less conservative and the approximate horizon may be closer to the optimal horizon corresponding to the true performance. Nonetheless, an upper bound itself may not suffice to help determine a good prediction horizon, and the addition of a performance \textit{lower bound} can be beneficial.

\begin{remark}[Consistency \& Rationality]
	The bounds given in \eqref{eq:performance_ce_mpc_affine} and \eqref{eq:performance_ce_mpc_linear} are both consistent in that they degenerate to 
	\begin{equation}
		\label{eq:performance_degenerate}
		J_{\infty}(x;\mu_{N,\hpara}, \tpara) \leq \frac{1}{1 - \epsilon_N}V_{\infty}(x;\tpara)
	\end{equation}
	when $\epara = 0$, recovering the performance bounds in \cite[Chapter 6]{grune2017nonlinear}, \cite{grune2008infinite, grune2012nmpc, kohler2023stability}. In other words, given that $\lim_{N \to +\infty}\epsilon_N = 0$, our performance bounds still advocate choosing a longer prediction horizon in the absence of model mismatch. On the other hand, the core factors (i.e., $\epsilon_N$, $\beta_{N}\big(\epara; x\big)$, $\beta^{\ast}_{N}(\epara)$, and $\alpha_N^\ast(\epara)$) that are used to build the final bounds only require knowing the mismatch level $\epara$ and other constants related to the properties of the cost functions (see Corollary \ref{coro:quadratic}) and $\hpara$, whereas $\tpara$ is never used. Therefore, our performance bounds are rational, i.e., they are compatible with the baseline motivation that the true model is inaccessible to the controller. It is worth noting that our analysis and results are also valid when the roles of $\tpara$ and $\hpara$ are swapped, i.e., when assuming $\tpara$ is known but $\hpara$ is unknown.  
\end{remark}
\begin{figure}[h]
\centering

\tikzset{every picture/.style={line width=0.75pt}} 
\resizebox{\linewidth}{!}{
\begin{tikzpicture}[x=0.75pt,y=0.75pt,yscale=-1,xscale=1]
	
	\draw  [color={rgb, 255:red, 245; green, 166; blue, 35 }  ,draw opacity=1 ][fill={rgb, 255:red, 248; green, 231; blue, 28 }  ,fill opacity=0.2 ] (0.33,1) -- (155.33,1) -- (155.33,36) -- (0.33,36) -- cycle ;
	\draw  [color={rgb, 255:red, 245; green, 166; blue, 35 }  ,draw opacity=1 ][fill={rgb, 255:red, 248; green, 231; blue, 28 }  ,fill opacity=0.2 ] (0.33,71) -- (155.33,71) -- (155.33,106) -- (0.33,106) -- cycle ;
	\draw  [color={rgb, 255:red, 245; green, 166; blue, 35 }  ,draw opacity=1 ][fill={rgb, 255:red, 248; green, 231; blue, 28 }  ,fill opacity=0.5 ] (0.33,141) -- (155.33,141) -- (155.33,176) -- (0.33,176) -- cycle ;
	\draw  [color={rgb, 255:red, 245; green, 166; blue, 35 }  ,draw opacity=1 ][fill={rgb, 255:red, 248; green, 231; blue, 28 }  ,fill opacity=0.2 ] (175.33,1) -- (330.33,1) -- (330.33,36) -- (175.33,36) -- cycle ;
	\draw  [color={rgb, 255:red, 245; green, 166; blue, 35 }  ,draw opacity=1 ][fill={rgb, 255:red, 248; green, 231; blue, 28 }  ,fill opacity=0.2 ] (175.33,71) -- (330.33,71) -- (330.33,106) -- (175.33,106) -- cycle ;
	\draw  [color={rgb, 255:red, 245; green, 166; blue, 35 }  ,draw opacity=1 ][fill={rgb, 255:red, 248; green, 231; blue, 28 }  ,fill opacity=0.5 ] (175.33,141) -- (330.33,141) -- (330.33,176) -- (175.33,176) -- cycle ;
	
	\draw [color={rgb, 255:red, 28; green, 128; blue, 128 }  ,draw opacity=1 ][fill={rgb, 255:red, 185; green, 185; blue, 185 }  ,fill opacity=1 ] [dash pattern={on 4.5pt off 4.5pt}]  (50.83,36) -- (50.83,70)(47.83,36) -- (47.83,70)(50.83,106) -- (50.83,132)(47.83,106) -- (47.83,132);
	\draw [shift={(49.33,140)}, rotate = 270] [fill={rgb, 255:red, 28; green, 128; blue, 128 }  ,fill opacity=1 ][line width=0.08]  [draw opacity=0] (10.2,-4.8) -- (0,0) -- (10.2,4.8) -- cycle    ;
	
	\draw [color={rgb, 255:red, 128; green, 128; blue, 128 }  ,draw opacity=1 ][fill={rgb, 255:red, 155; green, 155; blue, 155 }  ,fill opacity=1 ] (75.83,36) -- (75.83,62)(72.83,36) -- (72.83,62) ;
	\draw [shift={(74.33,70)}, rotate = 270] [fill={rgb, 255:red, 128; green, 128; blue, 128 }  ,fill opacity=1 ][line width=0.08]  [draw opacity=0] (10.2,-4.8) -- (0,0) -- (10.2,4.8) -- cycle    ;
	
	\draw [color={rgb, 255:red, 128; green, 128; blue, 128 }  ,draw opacity=1 ][fill={rgb, 255:red, 155; green, 155; blue, 155 }  ,fill opacity=1 ] (75.83,106) -- (75.83,132)(72.83,106) -- (72.83,132) ;
	\draw [shift={(74.33,140)}, rotate = 270] [fill={rgb, 255:red, 128; green, 128; blue, 128 }  ,fill opacity=1 ][line width=0.08]  [draw opacity=0] (10.2,-4.8) -- (0,0) -- (10.2,4.8) -- cycle    ;
	
	\draw [color={rgb, 255:red, 128; green, 128; blue, 128 }  ,draw opacity=1 ][fill={rgb, 255:red, 155; green, 155; blue, 155 }  ,fill opacity=1 ] (252.83,36) -- (252.83,62)(249.83,36) -- (249.83,62) ;
	\draw [shift={(251.33,70)}, rotate = 270] [fill={rgb, 255:red, 128; green, 128; blue, 128 }  ,fill opacity=1 ][line width=0.08]  [draw opacity=0] (10.2,-4.8) -- (0,0) -- (10.2,4.8) -- cycle    ;
	
	\draw [color={rgb, 255:red, 128; green, 128; blue, 128 }  ,draw opacity=1 ][fill={rgb, 255:red, 155; green, 155; blue, 155 }  ,fill opacity=1 ] (252.83,106) -- (252.83,132)(249.83,106) -- (249.83,132) ;
	\draw [shift={(251.33,140)}, rotate = 270] [fill={rgb, 255:red, 128; green, 128; blue, 128 }  ,fill opacity=1 ][line width=0.08]  [draw opacity=0] (10.2,-4.8) -- (0,0) -- (10.2,4.8) -- cycle    ;
	
	\draw [color={rgb, 255:red, 128; green, 128; blue, 128 }  ,draw opacity=1 ][fill={rgb, 255:red, 155; green, 155; blue, 155 }  ,fill opacity=1 ]   (155.33,155.5) -- (169.33,155.5)(155.33,158.5) -- (169.33,158.5) ;
	\draw [shift={(175.33,157)}, rotate = 180] [fill={rgb, 255:red, 128; green, 128; blue, 128 }  ,fill opacity=1 ][line width=0.08]  [draw opacity=0] (10.2,-4.8) -- (0,0) -- (10.2,4.8) -- cycle    ;
	
	\draw [color={rgb, 255:red, 128; green, 128; blue, 128 }  ,draw opacity=1 ][fill={rgb, 255:red, 155; green, 155; blue, 155 }  ,fill opacity=1 ]   (175.41,32.05) -- (142.99,65.32)(175.26,27.95) -- (140.84,63.22) ;
	\draw [shift={(136.33,70)}, rotate = 314.26] [fill={rgb, 255:red, 128; green, 128; blue, 128 }  ,fill opacity=1 ][line width=0.08]  [draw opacity=0] (10.2,-4.8) -- (0,0) -- (10.2,4.8) -- cycle    ;
	
	\draw (22,3) node [anchor=north west][inner sep=0.75pt]   [align=left] {\begin{minipage}[lt]{78.69pt}\setlength\topsep{0pt}
			\begin{center}
				\textbf{Nominal} Stability\\(Proposition \ref{prop:nominal_energy_decreasing})
			\end{center}
			
	\end{minipage}};
	\draw (175,3) node [anchor=north west][inner sep=0.75pt]   [align=left] {\begin{minipage}[lt]{115.49pt}\setlength\topsep{0pt}
			\begin{center}
				\textbf{Consistent} Error\\ Matching (Proposition \ref{prop:error_matching})
			\end{center}
			
	\end{minipage}};
	\draw (0,73) node [anchor=north west][inner sep=0.75pt]   [align=left] {\begin{minipage}[lt]{115.49pt}\setlength\topsep{0pt}
			\begin{center}
				\textbf{One-step-ahead State} \\ Perturbation (Proposition \ref{prop:parametric_state_perturbation})
			\end{center}
			
	\end{minipage}};
	\draw (0,144) node [anchor=north west][inner sep=0.75pt]   [align=left] {\begin{minipage}[lt]{115.78pt}\setlength\topsep{0pt}
			\begin{center}
				CE-MPC \textbf{Closed-loop} \\ \textbf{Stability} (Theorem \ref{thm:stability_ce_mpc})
			\end{center}
			
	\end{minipage}};
	\draw (175,73) node [anchor=north west][inner sep=0.75pt]   [align=left] {\begin{minipage}[lt]{115.74pt}\setlength\topsep{0pt}
			\begin{center}
				\textbf{Parameter} Perturbation \\ (Propositions \ref{prop:general_parameter_perturbation} and \ref{prop:linear_parameter_perturbation})
			\end{center}
			
	\end{minipage}};
	\draw (175,144) node [anchor=north west][inner sep=0.75pt]   [align=left] {\begin{minipage}[lt]{115.45pt}\setlength\topsep{0pt}
			\begin{center}
				CE-MPC \textbf{Closed-loop} \\ \textbf{Performance} (Theorem \ref{thm:affine_performance_ce_mpc})
			\end{center}
			
	\end{minipage}};

\end{tikzpicture}
}
\caption{Analysis pipeline for the stability and performance of certainty-equivalence model predictive control (CE-MPC).}
\vspace{-12pt}
\label{fig:pipeline}

\end{figure}

\section{Case Studies}
\label{sec:examples}
In this section, two case studies are provided to demonstrate the results of the previous sections as well as to gain insights from the proposed analysis framework. Particularly, we first provide a nonlinear system example with numerical simulation. Then, we specify our theoretical results for MPC of linear systems with parametric uncertainty and quadratic cost functions, serving as an important contribution in predictive control of uncertain linear systems. Compared to the previous results in the literature where affine and/or non-consistent bounds are provided \cite{moreno2022performance, liu2024stability}, or no constraint is considered \cite{lin2022bounded, lin2021perturbation}, we obtain the first consistent, rational, and competitive-ratio performance bound for input-constrained uncertain linear systems with multiplicative uncertainty. 
\subsection{Nonlinear System Simulation}
Consider the following second-order nonlinear system:
\begin{equation}
	\label{eq:nonlinear_numerical_example}
	\begin{aligned}
		& x_{k+1}[1] = -0.99x_{k}[2] \\
		& x_{k+1}[2] = w_1\tanh(x_{k}[1]) + w_2\tanh(x_{k}[2]) + bu_k
	\end{aligned},
\end{equation}
where $\tanh(\cdot)$ is the hyperbolic tangent function and $\para = [w_1, w_2, b]^\top$. This example is motivated by the popular neural network models (see \cite{pillonetto2025deep} and the references therein), where $\tanh(\cdot)$ is a widely adopted activation function\footnote{If the activation function is Lipschitz continuous, then the corresponding neural network model satisfies Assumption \ref{ass:sec_analysis:Lipschitz_dynamics}.}. The nominal parameter is $\hpara = [0.85, 0.995, 0.01]^\top$. The quadratic cost $\cost(x,u) = (x[1])^2 + (x[2])^2 + u^2$ is used, and $\cU = \{u \in \bR \mid -0.05 \leq u \leq 0.05\}$. It can be verified that $L_{f,x} = 0.995$ and $L_{f,u} = 0.01$ are valid nominal Lipschitz constants, and a linear mismatch function $\ldd(\mpara) = \mpara$ follows.

The simulation setting is given as follows. For a fixed $\epara$, multiple true models with different $\tpara$ satisfying $\|\tpara - \hpara\| = \epara$ are simulated, considering the worst case where $\delta(\para) = \varepsilon_{\para}$. Besides, the initial state will also be varied within the disk $(x[1])^2 + (x[2])^2 \leq 2$. The range of $N$ is restricted as $N \in \bI_{5:30}$, and the mismatch levels are $\epara = i \cdot 10^{-3}$ with $i \in \bI_{1:10}$ to ensure stability. As such, it holds that $\overline{L}_{f,x} = 1.005 > 1$ for $\epara = 0.01$. The nominal horizon is $N=10$. The MPC is solved using IPOPT with CasADi \cite{andersson2019} in Python 3.12.9\footnote{The code for the simulations is available at \url{https://github.com/lcrekko/CE-MPC}.}.

\begin{figure}[h]
	\centering
	\begin{subfigure}[b]{0.5\textwidth}
		\centering
		\includegraphics[width=\textwidth]{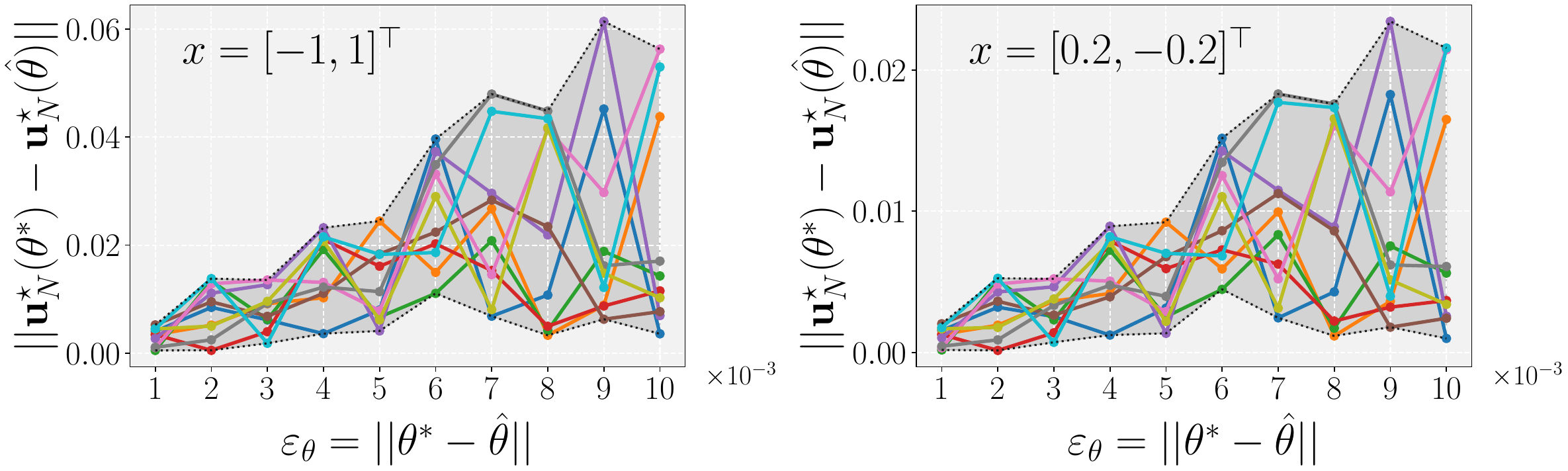}
		\caption{Input perturbation of typical $10$ scenarios for a fixed $\epara$.}
		\label{fig:input_pb_sub1}
	\end{subfigure}
	
	\vspace{1em} 
	
	\begin{subfigure}[b]{0.5\textwidth}
		\centering
		\includegraphics[width=\textwidth]{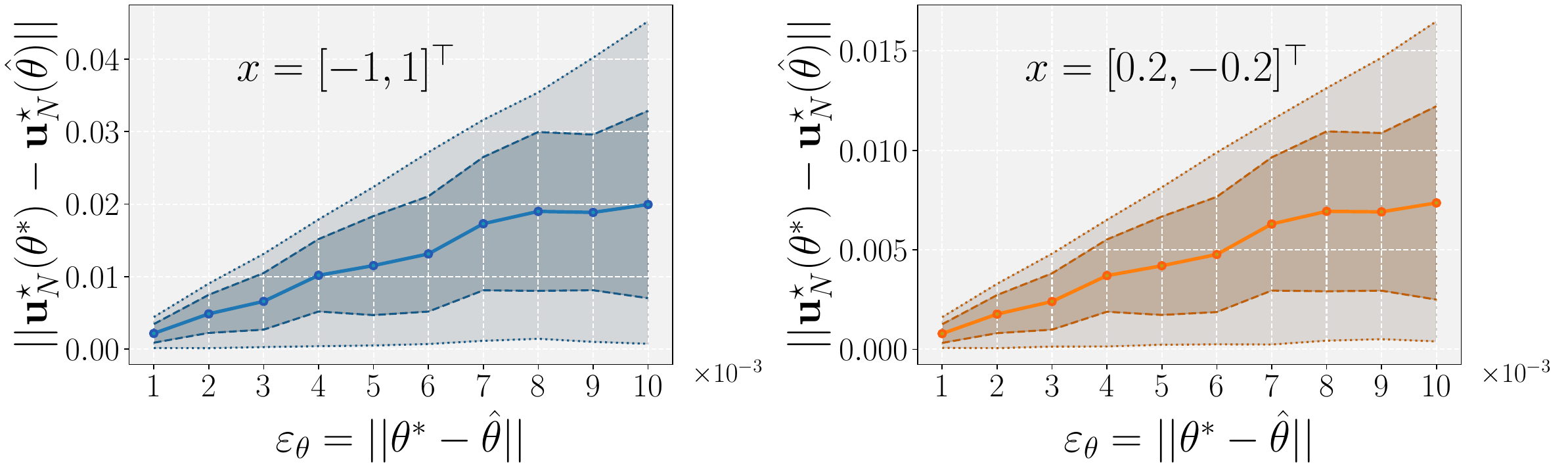}
		\caption{Input perturbation of $100$ scenarios for a fixed $\epara$. For the considered quantity $||\mathbf{u}^\star_N(\theta^\ast) - \mathbf{u}^\star_N(\hat{\theta})||$, the solid line is the mean value, the dashed lines are the mean plus (minus) variance, and the dotted lines are the max (min).}
		\label{fig:input_pb_sub2}
	\end{subfigure}
	
	\caption{Input perturbation simulation.}
	\label{fig:input_error_100}
\end{figure} 

\begin{figure}[h]
	\centering
	\includegraphics[width=0.6\linewidth]{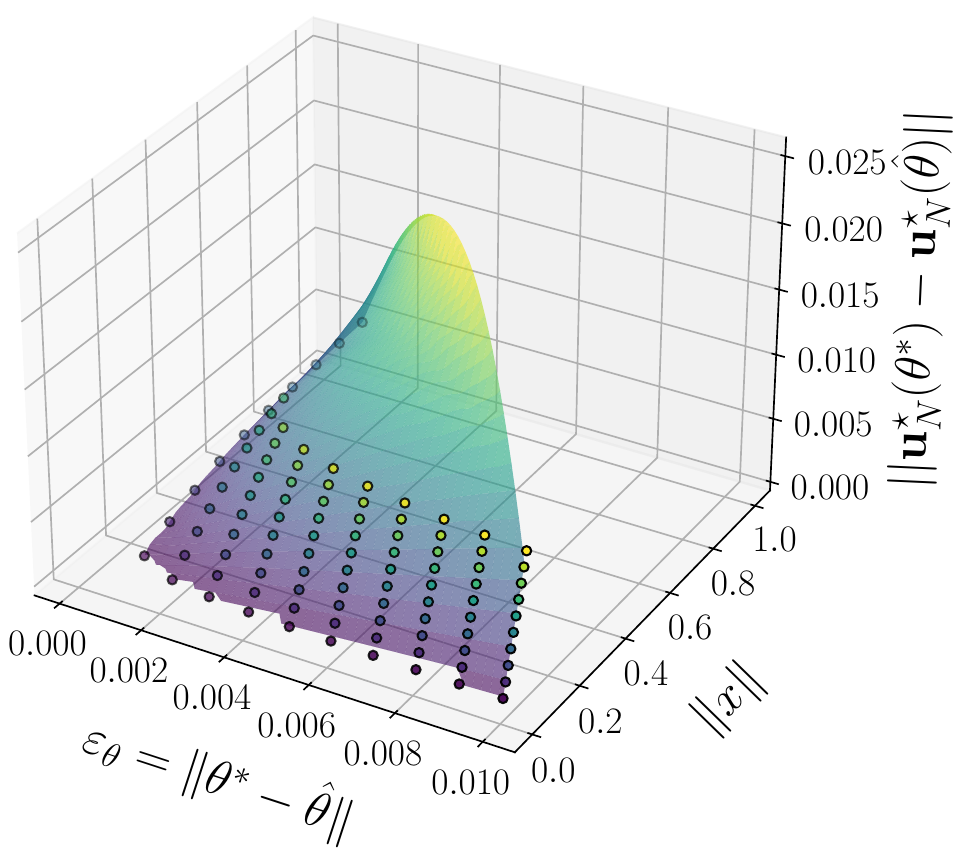}
	\caption{Scalable input perturbation of $100$ scenarios for a fixed $\epara$, where the maximum of $||\mathbf{u}^\star_N(\theta^\ast) - \mathbf{u}^\star_N(\hat{\theta})||$ is plotted.}
	\label{fig:scalable_input_error}
\end{figure}

The simulation results of input perturbation are given in Fig. \ref{fig:input_error_100} with $100$ scenarios of different $\tpara$ for each of the mismatch levels, from which a \textit{linear worst-case} bound with respect to $\epara$ can be observed, which is consistent with the linear bounds in Assumption \ref{ass:sec_analysis:ssosc_licq} and Lemma \ref{lm:sec_analysis:eds}. However, it should be noted that, though the bound of the perturbation is linear, the true perturbation is not necessarily linear and monotonously increasing (see Fig. \ref{fig:input_pb_sub1}).
Next, different initial conditions are considered such that multiple points on the level sets $r = \epara \|x\|$ are simulated with $r = i \cdot 10^{-3}$ with $i \in \bI_{1:10}$ such that a sufficient amount of points within the disk $(x[1])^2 + (x[2])^2 \leq \sqrt{2}$ are covered. The $3$-D visualization of the result is given in Fig. \ref{fig:scalable_input_error}, where a \textit{linear worst-case bound} with respect to the product $\epara \|x\|$ can be observed, validating Assumption \ref{ass:local_scalable_input_perturbation_bound} as well as Corollaries \ref{coro:scalable_eds} and \ref{coro:open_loop_parameter_perturbation}.

\begin{figure}[h]
	\centering
	\begin{subfigure}[b]{0.5\textwidth}
		\centering
		\includegraphics[width=\textwidth]{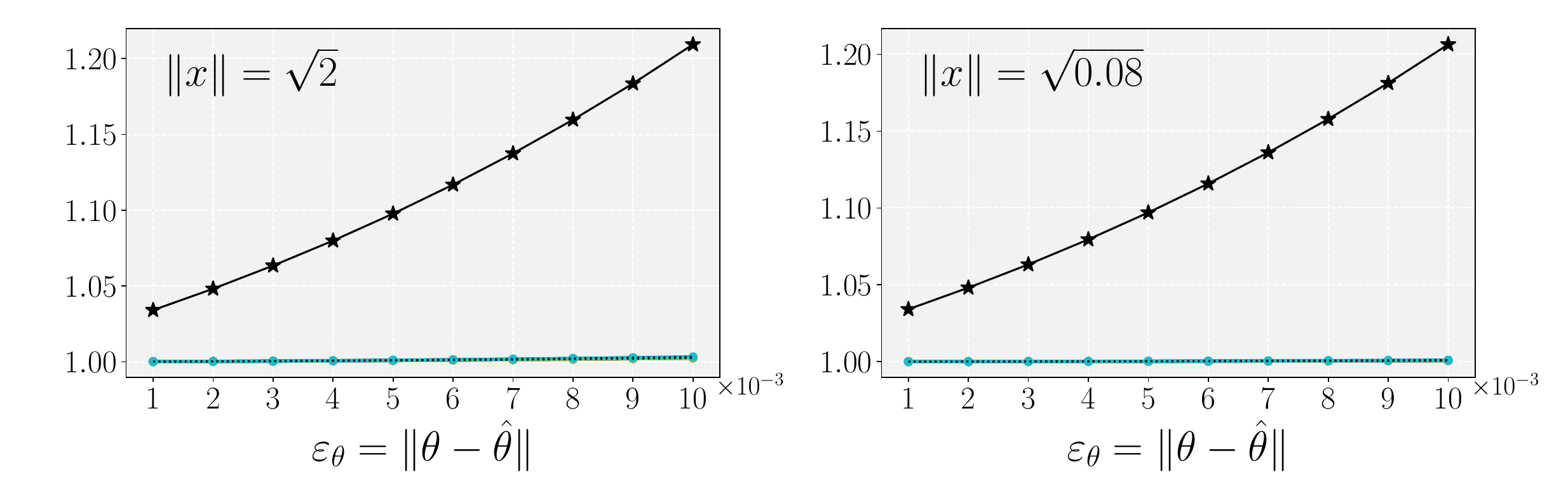}
		\caption{Numerical $\compratio{N}(\epara)$ (black line) and true competitive ratio of $100$ scenarios for a given $\epara$. $10$ different initial states are simulated for a given norm of the initial state.}
		\label{fig:cr_1}
	\end{subfigure}
	
	\vspace{1em} 
	
	\begin{subfigure}[b]{0.5\textwidth}
		\centering
		\includegraphics[width=\textwidth]{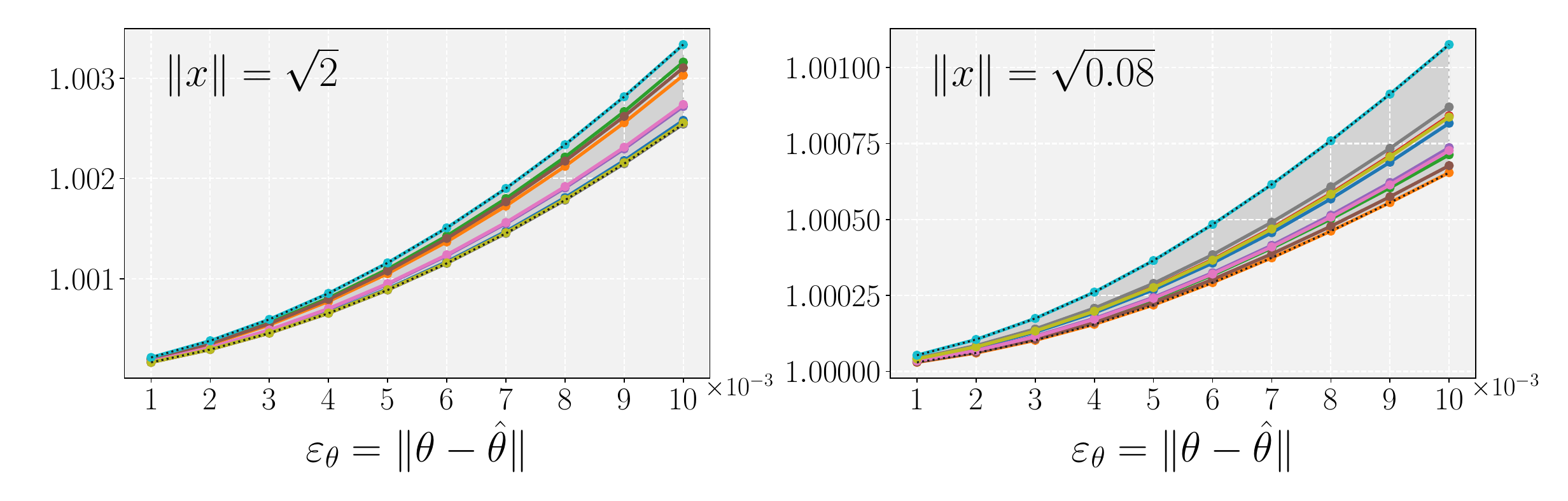}
		\caption{Detailed true competitive ratio curve of part (a) of this figure.}
		\label{fig:cr_2}
	\end{subfigure}
	
	\caption{Simulation of mismatch level $\epara$ variations.}
	\label{fig:competitive_ratio}
\end{figure}

Therefore, the competitive ratio bound \eqref{eq:performance_ce_mpc_linear} can be obtained for this nonlinear example, and Fig. \ref{fig:competitive_ratio} presents the bound $\compratio{N}(\epara)$ and the true competitive ratio. It is clear that the bound is valid but conservative. When $\epara$ is small, the true ratio grows linearly, and for relatively larger $\epara$ (approximately $\epara \geq 0.005$), the \textit{quadratic growth} is manifest, which coincides with the growth rate of the bound $\compratio{N}(\epara)$, thus supporting Propositions \ref{prop:parametric_state_perturbation} and \ref{prop:linear_parameter_perturbation} and the growth rate in Corollary \ref{coro:asymptotic_competitive_ratio}.

\begin{figure}[h]
	\centering
	\begin{subfigure}[b]{0.5\textwidth}
		\centering
		\includegraphics[width=\textwidth]{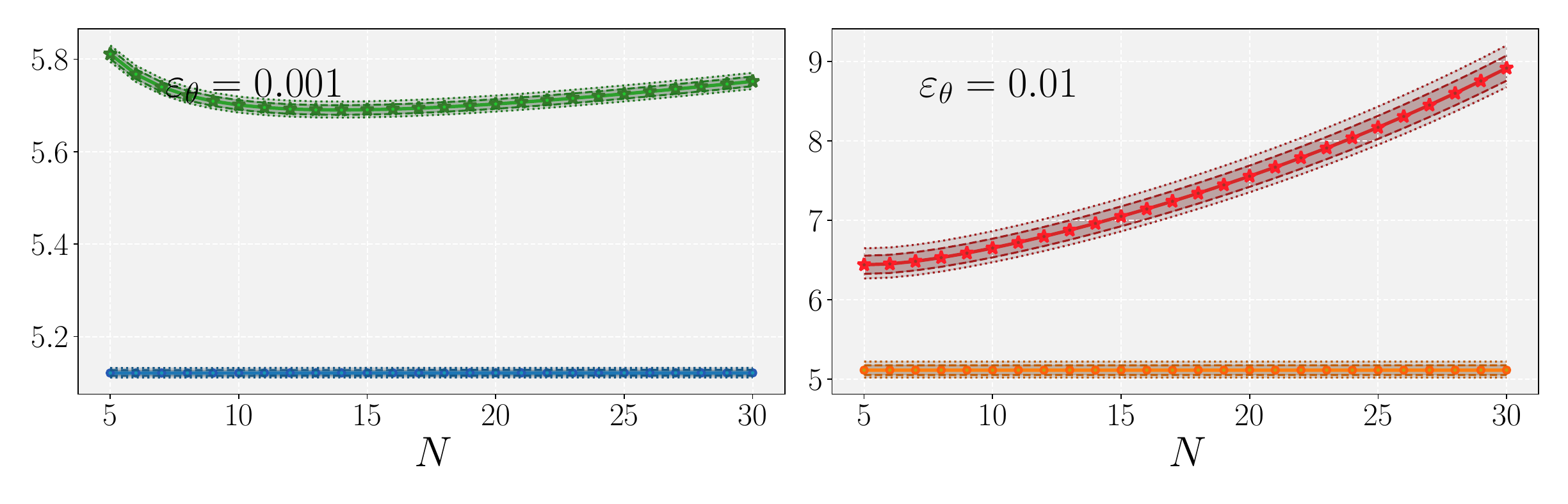}
		\caption{Performance bound (green and red curves) $\compratio{N}(\epara)V_{\infty}(x;\tpara)$ and true performance (blue and orange) $J_{\infty}(x;\mu_{N,\hat{\theta}},\theta^\ast)$ simulation of $100$ scenarios. The solid line is the mean, the dashed lines are the mean plus (minus) variance, and the dotted lines are the max (min).}
		\label{fig:horizon_1}
	\end{subfigure}
	
	\vspace{1em} 
	
	\begin{subfigure}[b]{0.5\textwidth}
		\centering
		\includegraphics[width=\textwidth]{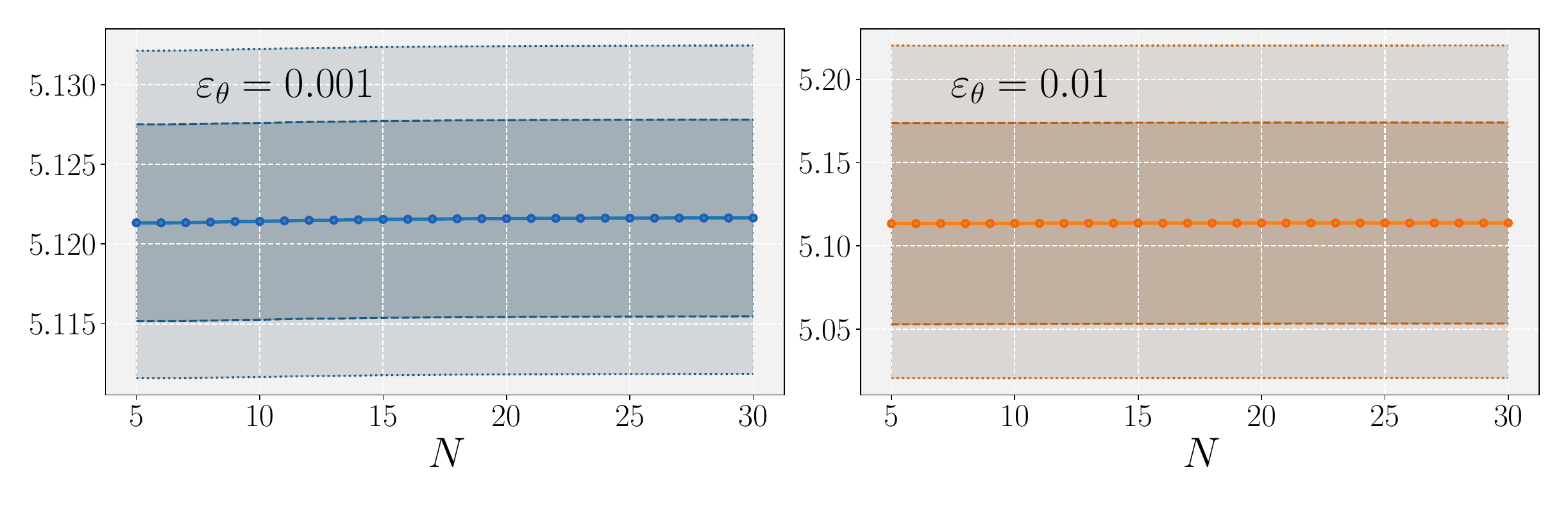}
		\caption{Detailed true performance curves of $100$ scenarios in part (a) of this figure.}
		\label{fig:horizon_2}
	\end{subfigure}
	
	\caption{Simulation of horizon $N$ variations}
	\label{fig:horizon_varation}
\end{figure}

The conservatism of \eqref{eq:stability_ce_mpc}, together with $\overline{L}_{f,x} > 1$, makes the simulation of the limit behavior when $N \to \infty$ difficult and meaningless since only infinitesimal mismatch can be considered (e.g., $\epara$ is of magnitude $10^{-9}$ or even smaller). Nonetheless, the considered range for $N$ is sufficient to convey the core messages. The effect of horizon variation is shown in Fig.~\ref{fig:horizon_varation} for different $\tpara$. It can be observed from Fig.~\ref{fig:horizon_1} that the derived performance bound in \eqref{eq:performance_ce_mpc_linear} is valid. However, as in Fig.~\ref{fig:horizon_2}, the true performance changes more significantly when varying $\tpara$ (i.e., the direction of the mismatch vector $\hpara - \tpara$) for a fixed $\epara$ compared to varying the horizon. Unfortunately, this vector sensitivity cannot be explained by our analysis since we focus on the worst-case bounds. In general, performance bounds do not reveal the true performance, which, especially for horizon variations, may not have significant changes as the bounds suggest (see also the examples in \cite{grune2008infinite, grune2012nmpc, grune2017nonlinear, kohler2021stability, kohler2023stability, giselsson2013feasibility}), and also may behave arbitrarily. Therefore, the bounds are only valuable in providing \textit{qualitative trends} on the worst-case performance. 


\vspace*{-0.3cm}
\subsection{Linear Quadratic Predictive Control}
Consider parametric linear systems of the form
\begin{equation}
	\label{eq:linear_system_model}
	x_{t+1} = A(\para)x_t + B(\para)u_t,
\end{equation}
where the matrices $A(\para)$ and $B(\para)$ are parameterized by $\para$. Without loss of generality, we assume there exist constants $\eab > 0$ such that for all $\para \in \setpara$, it holds that
\begin{equation}
	\label{eq:linear_matrices_basic_bound}
	\max\big\{\|A(\para) - A(\tpara)\|, \|B(\para) - B(\tpara)\|\big\} \leq \eab \mpara,
\end{equation}
which further implies that, for linear systems, the mismatch Lipschitz function is linear in $\mpara$, i.e, 
\begin{equation}
	\label{eq:mismatch_lipschitz}
	\ldd(\mpara) = \eab \mpara.
\end{equation}
It should be noted that the linear bound is valid for the widely adopted polytopic uncertainty \cite{fukushima2007adaptive, fleming2014robust, moreno2022performance, calafiore2012robust, bujarbaruah2021simple, arcari2023stochastic} and the more general normed uncertainty modeling \cite{shi2025suboptimality, liu2024stability, mania2019certainty} for linear systems with multiplicative uncertainty.
The quadratic cost function $\cost$ is chosen as
\begin{equation}
	\label{eq:cost_quadratic}
	\cost(x, u) = \nQ{x} + \nR{u},
\end{equation}
where $Q, R \succ 0$, and the optimized cost follows as $\cost^\star(x) = \nQ{x}$. The system described by the model in \eqref{eq:linear_system_model} satisfies Assumption \ref{ass:sec_analysis:Lipschitz_dynamics} with $L_{f,x} = \max_{\para \in \setpara} \|A(\para)\| =: \nA$ and $L_{f,u} = \max_{\para \in \setpara} \|B(\para)\| =: \nB$. The input constraint set is polytopic, i.e., $g_u(u) = E_u u$ with $E_u \in \bR^{c_u \times m}$. Besides, it is assumed that the pair $(A(\para), B(\para))$ is stabilizable for all $\para \in \setpara$. It is obvious that this linear quadratic control problem with polytopic input constraints satisfies Assumptions \ref{ass:sec_pre:basic_problem} and \ref{ass:sec_prob:terminal_cost}. Moreover, Assumption \ref{ass:local_scalable_input_perturbation_bound} is also valid (see Remark \ref{rmk:separable_local_input_perturbation}). Specifically, \eqref{eq:local_scalable_input_perturbation_bound} holds locally for the unconstrained linear quadratic control problem, i.e., $u^\star_{k|t}(x;\para) = K_{k,N}(\para)x$ (see Appendix \ref{appendix:E1}). In summary, a competitive ratio can be obtained for the linear quadratic control problem following \eqref{eq:bound_competitive_ratio}, as stated in the following corollary.

\begin{corollary}
	\label{coro:linear_quadratic_control}
	For the linear quadratic model predictive control problem specified in \eqref{eq:linear_system_model} to \eqref{eq:cost_quadratic}, let Assumptions \ref{ass:pre:initial_state_roa}, \ref{ass:sec_prob:nominal_cost_controllability}, \ref{ass:sec_analysis:ssosc_licq}, and the stability condition $\epsilon_N + \alpha_{\mathrm{LQ}}^\ast(\mpara) < 1$ hold for some $N \geq \underline{N}_0$ and $\epara \geq 0$. Then there exists $\omelq > 0$ such that $\sueds{k}^\ast$ in \eqref{eq:subsec_para_perturbation:refined_eds_separable} admits a linear upper bound $\sueds{k}^\ast(\mpara) \leq \omelq \mpara$ for all $k \in \bI_{0:N-1}$, and a competitive-ratio performance bound can be obtained, i.e.,
	\begin{equation}
		\label{eq:performance_ce_mpc_lq}
		J_{\infty}(x;\mu_{N,\hpara}, \tpara) \leq \compratiolq{N}(\epara)V_{\infty}(x;\tpara),
	\end{equation}
	with
	\begin{equation}
		\label{eq:bound_competitive_ratio_lq}
		\compratiolq{N}(\epara) = \frac{1 + \min\{\umpc, \beta_{\mathrm{LQ}}^\ast(\epara)\}}{1 - \big(\epsilon_N + \alpha_{\mathrm{LQ}}^\ast(\epara)\big)},
	\end{equation}
	where both $\alpha_{\mathrm{LQ}}^\ast(\epara)$ and $\beta_{\mathrm{LQ}}^\ast(\epara)$ are quadratic functions of the form $f(z) = az^2 + bz$ (see Propositions \ref{prop:parametric_state_perturbation} and \ref{prop:linear_parameter_perturbation}), and $\epsilon_N = \bigo{\rho^{N}}$ for some $\rho \in (0,1)$ \cite{liu2024stability}.
\end{corollary}


\section{Conclusions and future work}
\label{sec:conclusion} 
This paper has studied certainty-equivalence MPC for nonlinear systems with multiplicative parametric uncertainty and input constraints, focusing on the joint effect of model mismatch and the prediction horizon on the stability and infinite-horizon performance. Stability condition and suboptimality bounds have been derived. A key contribution is a perturbation analysis of the MPC value function that does not use Lipschitz cost, making the results suitable for the popular quadratic costs. The general results are further specified for linear quadratic control, serving as the first competitive-ratio performance bound for predictive control of uncertain linear systems with input constraints.

Future research directions include obtaining less conservative perturbation bounds, extensions to economic MPC and/or systems with additive disturbances and state constraints, investigating vector sensitivity when varying the direction of the unknown parameter under a fixed mismatch level, and deriving horizon-dependent performance lower bound for CE-MPC.

\section*{Acknowledgment}
The authors are grateful for the helpful discussions with Dr. Sungho Shin from MIT and Mr. Kanghui He from TU Delft.

\appendices
{\small
\section{Details of Section \ref{subsec_thereticalAnalysis:preparatory}}
\label{appendix:A--SecIV-A}

\subsection{Proof of Lemma \ref{lm:dynamic_difference}}
\label{appendix:A0}
\begin{proof}
	$\Delta_f(\cdot,\cdot;\para)$ is Lipschitz continuous since $f$ is Lipschitz continuous. Thus, noting $\Delta_f(x,u;\hpara) = \mathbf{0}$, there exist positive definite functions $\lddx', \lddu'$ such that $\|\Delta_f(x',u';\para) - \Delta_f(x'',u'';\para)\| \leq \lddx'(\mpara)\|x' - x''\| + \lddu'(\mpara)\|u' - u''\|$, $\forall x',x''\in \cX, u',u'' \in \cU$, and both $\lddx'$ and $\lddu'$ take finite values. Besides, it holds that $\Delta_f(\mathbf{0},\mathbf{0};\para) = \mathbf{0}$, then $\|\Delta_f(x',u';\para)\| \leq \lddx'(\mpara)\|x'\| + \lddu'(\mpara)\|u'\|$, noting $x''=\bzero$ and $u''=\bzero$. Given $\epsilon > 0$, construct $\lddx$ and $\lddu$, respectively, as $\lddx(\delta) = \max_{0\leq\delta'\leq\delta}\lddx'(\delta')+\epsilon \delta$ and $\lddu(\delta) = \max_{0\leq\delta'\leq\delta}\lddu'(\delta')+\epsilon \delta$. Then $\lddx$ and $\lddu$ are strictly increasing, and they satisfy $\lddx(\delta) \geq \lddx'(\delta)$ and $\lddu(\delta) \geq \lddu'(\delta)$, respectively, for all $\delta \in [0,\epara]$, which completes the proof.
\end{proof}

\subsection{Proof of Proposition \ref{prop:error_matching}}
\label{appendix:A2}
\begin{proof}
The following proof is inspired by \cite[Lemma 13]{liu2024stability}. First, when $\lddx(\mpara), \lddu(\mpara) > 0$,
\begin{align}
	\label{eq:A2_error_lb}
	& \|\dsys\|^2 \notag \\
	\overset{\eqref{eq:difference_lipschitz_dynamics}}{\leq} & \left(\lddx(\mpara)\|x\| + \lddu(\mpara)\|u\|\right)^2 \notag \\
	\leq & \left(1 + \frac{\lddu^2(\mpara)\mlx}{\lddx^2(\mpara)\mlu}\right)\lddx^2(\mpara)\|x\|^2 + \notag \\ 
	& \hspace*{5ex} \left(1 + \frac{\lddx^2(\mpara)\mlu}{\lddu^2(\mpara)\mlx}\right)\lddu^2(\mpara)\|u\|^2 \notag \\
	\overset{\eqref{eq:quadratic_lower}}{\leq} &  \; 2\left(\frac{\lddx^2(\mpara)}{\mlx} + \frac{\lddu^2(\mpara)}{\mlu}\right)[\ell_x(x) + \ell_u(u)], \notag \\
	\overset{\eqref{eq:big_O_ease_difference_upper_bound}}{\leq} & \; \underbrace{2\left(\frac{1}{\mlx} + \frac{1}{\mlu}\right)}_{:=\mec}\ldd^2(\mpara)\ell(x,u),
\end{align}
where the second inequality is by Cauchy-Schwarz and Young's inequalities. If $\lddx(\mpara) = \lddu(\mpara) = 0$, then $\|\dsys\| = 0$, and \eqref{eq:A2_error_lb} holds trivially.
\end{proof}

\section{Details of Section \ref{subsec_theoreticalAnalysis:state_perturbation_mpc}}
\label{appendix:B--SecIV-B}

\subsection{Proof of Proposition \ref{prop:nominal_energy_decreasing}}
\label{appendix:B1}
\begin{proof}
	The property \eqref{eq:nominal_bound} can be simply verified by definition and Assumption \ref{ass:sec_prob:nominal_cost_controllability}, and \eqref{eq:nominal_decreasing} can be established using the case distinction method \cite{kohler2021stability, kohler2023stability}, where our distinction condition is 
	\begin{equation}
		\ell^\star(\xi^\star_N(x;\hpara)) \leq \frac{(1+\umpc)\mpc_{N}}{(N-1)\epsilon_F + N + \umpc}\ell^\star(x).
	\end{equation}
	Since the reasoning follows that of \cite[Appendix A-(1)-Part I]{kohler2023stability}, we omit the remaining details here for brevity, and existing proofs in the literature \cite{kohler2023stability}, \cite[Chapter 4]{grune2017nonlinear} are sufficient to support our proof. Next, we give the proof of \eqref{eq:next_state_value_upper_bound} as follows:
	\begin{align}
		\label{eq:B3_energy_decreasing}
		\;V_N(\xp{\hpara};\hpara)
		\overset{\eqref{eq:nominal_decreasing}}{\leq} & \;V_N(x;\hpara) - (1 - \epsilon_N)\cost(x, \mu_{N,\hpara}(x))\notag  \\ 
		\overset{\eqref{eq:sec_prob:mpc_linear_bound}}{\leq} &  \;(1 + \mpc_N)\cost^\star(x) - (1 - \epsilon_N)\cost(x, \mu_{N,\hpara}(x))\notag  \\ 
		 \overset{\eqref{eq:sec_pre:optimized_stage_cost}}{\leq}&  \;(\mpc_N + \epsilon_N)\cost(x, \mu_{N,\hpara}(x)),
	\end{align}
	which completes the proof. 
\end{proof}

\subsection{Proof of Lemma \ref{lm:open_loop_state_perturbation}}
\label{appendix:B2}
\begin{proof}
	The proof employs induction, and we omit the details for brevity. When selecting $x' = \xp{\hpara}$ and $x'' = \xp{\para}$, a slightly modified version of \eqref{eq:lm4_open_loop_state_perturbation_bound} is
	\begin{equation}
		\label{eq:B2_final_conclusion}
		\|\xi^\star_{k}(\xp{\hpara};\hpara) - \xi^\star_{k}(\xp{\para};\hpara)\| \leq \Gamma_{k}\|\xp{\hpara} - \xp{\para}\|,
	\end{equation}
	which will be used in the proof of Proposition \ref{prop:parametric_state_perturbation}.
\end{proof}

\subsection{Proof of Proposition \ref{prop:parametric_state_perturbation}}
\label{appendix:B3}
\begin{proof}
	For notational convenience, we define the shorthand notations $\Delta \xi^+_k := \xi^\star_{k}(\xp{\hpara};\hpara) - \xi^\star_{k}(\xp{\para};\hpara)$ and $\Delta u^+_k := u^\star_{k}(\xp{\hpara};\hpara) - u^\star_{k}(\xp{\para};\hpara)$.
	Leveraging \eqref{eq:cost_perturbation}, it holds that
	\begin{multline}
		\label{eq:B3_initial}
		|V_{N}(\xp{\para};\hpara) - V_{N}(\xp{\hpara};\hpara)| \leq \frac{\llx}{2}\sum^N_{k=0}\dxip{k}^2 + \\ \frac{\llu}{2}\sum^{N-1}_{k=0}\dup{k}^2	+ 
		\sum^N_{k=0}\llx \|\xiph{k}\|\dxip{k} + \\ \sum^{N-1}_{k=0}\llu \|\uph{k}\|\dup{k}.
	\end{multline}
	First, the terms related to $\dxip{k}^2$ and $\dup{k}^2$ are simplified as
	\begin{align}
		\label{eq:B3_second_order}
		& \hspace*{3ex} \frac{\llx}{2}\sum^N_{k=0}\dxip{k}^2 + \frac{\llu}{2}\sum^{N-1}_{k=0}\dup{k}^2 \notag \\
		& \overset{\eqref{eq:B2_final_conclusion},\eqref{eq:subsec_eds:state_sensitivity_input}}{\leq} \left[\frac{\llx}{2}\sum^{N}_{k=0}\Gamma^2_k + \frac{\llu \edsc^2}{2(1 - \edsr^2)}\right]\|\xp{\hpara} - \xp{\para}\|^2 \notag \\
		& \hspace*{-1ex} \overset{\eqref{eq:error_consistency}}{\leq}  \underbrace{\mec\left[\frac{\llx}{2}\sum^{N}_{k=0}\Gamma^2_k \splus \frac{\llu \edsc^2}{2(1 - \edsr^2)}\right]}_{:=\pi_{\alpha,N,2}} \ldd^2(\mpara)\cost(x,\law). \hspace*{-1ex}
	\end{align}
	On the other hand, the terms associated with $\dxip{k}$ and $\dup{k}$ can be further written as
	\begin{align}
		\label{eq:B3_first_order}
		& \hspace*{3ex} \sum^N_{k=0}\llx \|\xiph{k}\|\dxip{k} + \notag \\ 
		& \hspace*{13ex} \sum^{N-1}_{k=0}\llu \|\uph{k}\|\dup{k} \notag \\
		& \leq \left\{ \left[\sum^N_{k=0}\mlx\|\xiph{k}\|^2 \splus \sum^{N-1}_{k=0}\mlu\|\uph{k}\|^2 \right]\cdot \right. \notag \\
		& \hspace*{8ex} \left. \left[\sum^{N}_{k=0}\frac{\llx^2}{\mlx}\dxip{k}^2 + \sum^{N-1}_{k=0}\frac{\llu^2}{\mlu}\dup{k}^2 \right] \right\}^{\frac{1}{2}} \notag \\
		& \overset{\eqref{eq:subsec_eds:state_sensitivity_input},\eqref{eq:cost_perturbation},\eqref{eq:B2_final_conclusion}}{\leq} \hspace*{-1ex} \left\{ 2\left[\sum^N_{k=0}\cost_x(\xiph{k}) \splus \sum^{N-1}_{k=0}\cost_u(\uph{k}) \right]\hspace*{-0.5ex} \cdot \right. \notag \\
		& \hspace*{5ex} \left. \left(\frac{\llx^2}{\mlx}\sum^{N}_{k=0}\Gamma^2_k + \frac{\llu^2\edsc}{\mlu(1 - \edsr^2)}\right)\|\xp{\hpara} - \xp{\para}\|^2 \right\}^{\frac{1}{2}} \notag \\
		& \overset{\eqref{eq:error_consistency},\eqref{eq:next_state_value_upper_bound}}{\leq} \underbrace{\left[2\mec(\umpc + \epsilon_N)\left(\frac{\llx^2}{\mlx}\sum^{N}_{k=0}\Gamma^2_k + \frac{\llu^2\edsc}{\mlu(1 - \edsr^2)}\right)\right]^{\frac{1}{2}}}_{:=\pi_{\alpha,N,1}}\cdot \notag \\
		& \hspace*{30ex} \ldd(\mpara)\cost(x,\law),
	\end{align}
	where the second inequality is due to Cauchy-Schwartz inequality. Finally, substituting \eqref{eq:B3_second_order} and \eqref{eq:B3_first_order} into \eqref{eq:B3_initial} yields the final perturbation bound \eqref{eq:parametric_state_perturbation} with $\alpha^\ast_N$ given in \eqref{eq:alpha_detail_in_the_theorem}. 
\end{proof}

\section{Details of Section \ref{subsec_theoreticalAnalysis:parameter_perturbation_mpc}}
\label{appendix:C--SecIV-C}

\subsection{Proof of Lemma \ref{lm:scalable_eds}}
\label{appendix:C1}
\begin{proof}
	According to Assumption \ref{ass:sec_pre:basic_problem}, for any $\para \in \setpara$, it holds that $u^\star_{k}(\mathbf{0};\para) = 0$, which, based on \eqref{eq:subsec_eds:state_sensitivity_input}, implies $\|u^\star_{k}(x;\para)\| \leq \edsc\edsr^k\|x\|$ on $\mathcal{B}(\mathbf{0},R(\mathbf{0};\epara))$ . Thus, by the triangle inequality,
	\begin{equation}
		\label{eq:C1_local_state_bound}
		\|u^\star_{k}(x;\para) \sminus u^\star_{k}(x;\hpara)\| \leq 2\edsc\edsr^k\|x\|, \forall x \in \mathcal{B}(\mathbf{0},R(\mathbf{0};\epara)),
	\end{equation}
	which, after leveraging \eqref{eq:subsec_eds:parametric_sensitivity}, leads to
	\begin{multline}
		\label{eq:C1_local_combined_bound}
		\|u^\star_{k}(x;\para) - u^\star_{k}(x;\hpara)\| \leq \min\{\Lambda_N(k)\mpara,2\edsr^k\|x\|\}\edsc \leq \\ 2\edsc\min\{\frac{\mpara}{1-\edsr},\|x\|\}, \forall x \in \mathcal{B}(\mathbf{0},R(\mathbf{0};\epara)).
	\end{multline}
	In addition, for $x \in \xroa \setminus \mathcal{B}(\mathbf{0},R(\mathbf{0};\epara))$, \eqref{eq:subsec_eds:parametric_sensitivity} directly yields
	\begin{equation}
		\label{eq:C1_not_local}
		\|u^\star_{k}(x;\para) - u^\star_{k}(x;\hpara)\| \leq
	\frac{2\edsc}{(1-\edsr)R(\mathbf{0};\epara)}\mpara \|x\|.
	\end{equation}
	Combining \eqref{eq:C1_local_combined_bound} and \eqref{eq:C1_not_local} leads to the desired result \eqref{eq:subsec_para_perturbation:refined_eds}.
\end{proof}

\subsection{Proof of Lemma \ref{lm:open_loop_parameter_perturbation}}
\label{appendix:C2}
\begin{proof}
	For convenience, $\forall k \in \bI_{0:N}$, we define shorthand notations $\Delta \xi^\star_{k}(x) := \xi^\star_{k}(x;\para) - \xi^\star_{k}(x;\hpara)$, $\Delta u^\star_{k}(x) :=  u^\star_{k}(x;\para) - u^\star_{k}(x;\hpara)$, and $\Delta^\star_{f,k}(x;\hpara) := \Delta_f(\xi^\star_{k}(x;\hpara), u^\star_{k}(x;\hpara); \para)$. First, we establish the recursive relation between $\difxi{k+1}$ and $\difxi{k}$ as follows:
	\begin{align}
		\label{eq:C2_recursive_relation}
		& \;\difxi{k+1} \notag \\ 
		= & \;\|f(\xi^{\star}_k(x;\para), u^{\star}_k(x;\para); \para) \sminus f(\xi^{\star}_k(x;\hpara), u^{\star}_k(x;\hpara); \hpara)\| \notag \\
		\leq &\; \|f(\xi^{\star}_k(x;\para), u^{\star}_k(x;\para); \para) \sminus f(\xi^{\star}_k(x;\hpara), u^{\star}_k(x;\hpara); \para)\| + \notag \\
		&\; \hspace*{35ex} \dpairf{k} \notag \\
		& \hspace*{-2ex} \overset{\eqref{eq:lipschitz_dynamics}}{\leq}  \; L_{f,x}\difxi{k} + L_{f,u}\difu{k} + \dpairf{k},
	\end{align}
	where the second inequality is due to the triangle inequality. The result \eqref{eq:open_loop_parameter_perturbation_bound} is proven by induction, which is omitted for brevity.
\end{proof}

\subsection{Proof of Proposition \ref{prop:general_parameter_perturbation}}
\label{appendix:C3}
\begin{proof}
	Lemma \ref{lm:scalable_eds} provides an uniform input perturbation bound $\|u^\star_{k}(x;\hpara) - u^\star_{k}(x;\para)\| \leq \bo(\mpara; \|x\|)$. For convenience, in the following proof, $\bo$ will be used as a shorthand for $\bo(\mpara; \|x\|)$. Besides, the shorthand notations in Appendix \ref{appendix:C2} will be adopted as well. By \eqref{eq:cost_perturbation} and noting that $\difxi{0} = 0$, it holds that
	\begin{multline}
		\label{eq:C3_beginning}
		|V_N(x;\hpara) - V_N(x;\para)| \leq \frac{\llx}{2}\sum^N_{k=1}\difxi{k}^2 + \\ \frac{\llu}{2}\sum^N_{k=1}\difu{k}^2 + 
		\llx\sum^N_{{k=1}}\|\xin{k}\|\difxi{k} + \\ \llx\sum^N_{{k=1}}\|\uun{k}\|\difu{k},
	\end{multline}
	where the terms coupling with $\xin{k}$ and $\uun{k}$ can be further relaxed using the Cauchy-Schwartz inequality as
	\begin{multline}
		\label{eq:C3_crossterms}
		\hspace*{-2ex}\llx\sum^N_{{k=1}}\|\xin{k}\|\difxi{k} + \llx\sum^N_{{k=1}}\|\uun{k}\|\difu{k} \\
		\leq \left\{ \left[\sum^N_{k=1}\mlx\|\xin{k}\|^2 + \sum^{N-1}_{k=0}\mlu\|\uun{k}\|^2\right]\cdot\right. \\
		\left. \left[\frac{\llx^2}{\mlx}\sum^N_{k=1}\difxi{k}^2 + \frac{\llu^2}{\mlu}\sum^{N-1}_{k=0}\difu{k}^2\right]\right\}^{\frac{1}{2}} \\
		\hspace*{-5ex}\overset{\eqref{eq:cost_quadratic}}{\leq}\sqrt{2}\left\{ \left[\sum^N_{k=1}\cost_x(\xin{k}) + \sum^{N-1}_{k=0}\cost_u(\uun{k})\right]\cdot\right. \\
		\left. \left[\frac{\llx^2}{\mlx}\sum^N_{k=1}\difxi{k}^2 + \frac{\llu^2}{\mlu}\sum^{N-1}_{k=0}\difu{k}^2\right]\right\}^{\frac{1}{2}}  \\
		\hspace*{-25ex}\leq \sqrt{2}\Bigg\{ \left[V_N(x;\hpara) - \cost^\star(x)\right]\cdot \\
		\left. \left[\frac{\llx^2}{\mlx}\sum^N_{k=1}\difxi{k}^2 + \frac{\llu^2}{\mlu}\sum^{N-1}_{k=0}\difu{k}^2\right]\right\}^{\frac{1}{2}} \\
		\hspace*{-40ex}\overset{\eqref{eq:sec_prob:mpc_linear_bound}}{\leq} \sqrt{2\umpc}\Bigg\{\cost^\star(x) \\
		\left. \left[\frac{\llx^2}{\mlx}\sum^N_{k=1}\difxi{k}^2 + \frac{\llu^2}{\mlu}\sum^{N-1}_{k=0}\difu{k}^2\right]\right\}^{\frac{1}{2}}.
	\end{multline}
	It is obvious that $\sum^{N-1}_{k=0}\difu{k}^2 \leq N\bo^2$. Besides, define $P(\blfx) = \sum^N_{k=1}\blfu^2\left(\sum^{k-1}_{i=0}\blfx^{i}\right)^2 + \sum^N_{k=1}\sum^{k-1}_{i=0}\blfx^{2i}$; then by the Cauchy-Schwartz inequality and \eqref{eq:error_consistency}, it can be obtained that
	\begin{equation}
		\label{eq:C3_state_summation}
		\hspace*{-1.3ex}\sum^N_{k=1}\difxi{k}^2 \leq P(\blfx)\left[\bo^2 \splus \mec \ldd^2(\mpara)\umpc\cost^\star(x)\right].\hspace*{-0.7ex}
	\end{equation}
	Substituting \eqref{eq:C3_crossterms} and \eqref{eq:C3_state_summation} into \eqref{eq:C3_beginning} leads to
	\begin{multline}
		\label{eq:C3_final}
		|V_N(x;\hpara) - V_N(x;\para)| \leq \underbrace{\left[\frac{\llx}{2}\umpc P(\blfx)\mec\cost^\star(x)\right]}_{:=\pi_{\beta,N,2}(x)}\ldd^2(\mpara) + \\
		\underbrace{\left[\frac{2\mec P(\blfx)}{\mlx}\right]^{\frac{1}{2}}\llx\umpc\cost^\star(x)}_{:=\pi_{\beta,N,1}(x)}\ldd(\mpara) + \\ \underbrace{\frac{1}{2}[\llx P(\blfx) + \llu N]}_{:=\zeta_{\beta,N,2}}\bo^2(\mpara; \|x\|) + \\
		\underbrace{[2\umpc\cost^\star(x)]^{\frac{1}{2}}[\frac{\llx}{\sqrt{\mlx}}P^{\frac{1}{2}}(\blfx) + \frac{\llx}{\sqrt{\mlx}}N^{\frac{1}{2}}]}_{:=\zeta_{\beta,N,1}(x)}\bo(\mpara; \|x\|).
	\end{multline}
	In addition, it can be verified that
	\begin{equation}
		\label{eq:C3_bigoP}
		P(\blfx) =
		\begin{cases}
			\bigo{N} & \blfx < 1 \\
			\bigo{N^3} & \blfx = 1 \\
 			\bigo{\blfx^{2N}} & \blfx > 1
		\end{cases},
	\end{equation}
	when $N \to \infty$, which, together with \eqref{eq:C3_final}, yields the asymptotic behavior given in Proposition \ref{prop:general_parameter_perturbation}.
\end{proof}

\subsection{Proof of Proposition \ref{prop:linear_parameter_perturbation}}
\label{appendix:C4}
Before presenting the proof, we first give additional comments regarding Corollaries \ref{coro:scalable_eds} and \ref{coro:open_loop_parameter_perturbation}. First, Corollary \ref{coro:scalable_eds} is a direct modification of Lemma \ref{lm:scalable_eds} based on Assumption \ref{ass:local_scalable_input_perturbation_bound}, and it can be proven using the same arguments as in Lemma \ref{lm:scalable_eds}. Then, Corollary \ref{coro:open_loop_parameter_perturbation} is derived from Lemma \ref{lm:open_loop_parameter_perturbation} based on Corollary \ref{coro:scalable_eds}, and its proof is almost similar to that of Lemma \ref{lm:open_loop_parameter_perturbation} given in Appendix \ref{appendix:C2} with the only difference being that the recursive relation \eqref{eq:C2_recursive_relation} is now modified as
	\begin{align}
	\label{eq:C4_recursive_relation}
	& \;\difxi{k+1} \notag \\ 
	= & \;\|f(\xi^{\star}_k(x;\para), u^{\star}_k(x;\para); \para) \sminus f(\xi^{\star}_k(x;\hpara), u^{\star}_k(x;\hpara); \hpara)\| \notag \\
	\leq &\; \dpaf{k} + \notag \\
	&\; \|f(\xi^{\star}_k(x;\para), u^{\star}_k(x;\para); \hpara) \sminus f(\xi^{\star}_k(x;\hpara), u^{\star}_k(x;\hpara); \hpara)\| \notag \\
	& \hspace*{-2ex} \overset{\eqref{eq:lipschitz_dynamics}}{\leq}  \; L_{f,x}\difxi{k} + L_{f,u}\difu{k} + \dpaf{k},
\end{align}
where $\Delta^\star_{f,k}(x;\para) := \Delta_f(\xi^\star_{k}(x;\para), u^\star_{k}(x;\para); \para),\;\forall k \in \bI_{0:N-1}$.
\begin{proof}
	The proof follows the same reasoning as that of Proposition \ref{prop:general_parameter_perturbation}, with the only difference being that the first decomposition (cf. \eqref{eq:C3_beginning}) should be constructed based on $\xit{k}$ and $\uut{k}$. The details are omitted for brevity and the final result is
	\begin{multline}
		\label{eq:C4_final}
		\hspace*{-2ex}|V_N(x;\hpara) - V_N(x;\para)| \leq \Bigg\{ \underbrace{\left[\frac{\llx}{2}\umpc P(\blfx)\mec\cost^\star(x)\right]}_{:=\pi^\ast_{\beta,N,2}}\ldd^2(\mpara) + \\
		\underbrace{\left[\frac{2\mec P(\blfx)\umpc}{\mlx}\right]^{\frac{1}{2}}\llx}_{:=\pi^\ast_{\beta,N,1}}\ldd(\mpara) + \\ \underbrace{(\overline{\eta}^\ast)^2\left[\frac{\llx}{\mlx} P(\blfx) + \frac{\llu}{\mlx} N\right]}_{:=\zeta^\ast_{\beta,N,2}}\delta^2(\para) + \\
		\underbrace{2\overline{\eta}^\ast\left[\frac{\llx}{\mlx}P^{\frac{1}{2}}(\blfx) + \frac{\llu}{\sqrt{\mlx\mlu}}N^{\frac{1}{2}}\right]}_{:=\zeta^\ast_{\beta,N,1}}\mpara\Bigg\} V_N(x;\para).
	\end{multline}
	The asymptotic behavior can be obtained similarly using \eqref{eq:C3_bigoP}.
\end{proof}

\section{Details in Section \ref{sec:theoreticalAnalysis:stability_suboptimality}}
\label{appendix:D}

\subsection{Proof of Theorem \ref{thm:stability_ce_mpc}}
\label{appendix:D1}
\begin{proof}
	The proof follows the standard RDP framework \cite{grune2017nonlinear, kohler2023stability}. Specifically, $V_N(\xp{\tpara}; \hpara) - V_N(x; \hpara)$ will be analyzed.
	\begin{align}
		\label{eq:D1_main_argument}
		& \hspace*{4ex} V_N(\xp{\tpara}; \hpara) - V_N(x; \hpara) \notag \\
		& \leq \; |V_N(\xp{\tpara}; \hpara) - V_N(\xp{\hpara}; \hpara)| + V_N(\xp{\hpara}; \hpara) - V_N(x; \hpara) \notag \\
		& \overset{\eqref{eq:nominal_decreasing}, \eqref{eq:parametric_state_perturbation}}{\leq}  -\big(1 - \alpha^\ast_N(\epara) - \epsilon_N\big)\cost(x, \mu_{N,\hpara}(x)),
	\end{align}
	Eq.~\eqref{eq:D1_main_argument} implies that $V_N(\xp{\tpara}; \hpara) - V_N(x; \hpara) \leq -\big(1 - \alpha^\ast_N(\epara) - \epsilon_N\big)\cost^\star(x)$, which, together with \eqref{eq:stability_ce_mpc} and \eqref{eq:nominal_bound}, proves that $V_N(\cdot;\hpara)$ is a Lyapunov function for the closed-loop system \eqref{eq:sec_prob:closed_loop_system}. In addition, performing the telescopic sum $T$ times on \eqref{eq:D1_main_argument} yields
	\begin{multline}
		\label{eq:D1_telescopic_sum}
		\hspace*{-2ex} \big(1 - \alpha^\ast_N(\epara) - \epsilon_N\big)\sum^{T-1}_{t=0}\cost(\ctx{\tpara}(t;x,\mu_{N,\hpara}), \ctu{\tpara}(t;x,\mu_{N,\hpara})) \\ 
		\leq V_N(x; \hpara) - V_N(\ctx{\tpara}(T;x,\mu_{N,\hpara}); \hpara) \leq V_N(x; \hpara).
	\end{multline}
	Take $T \to +\infty$, \eqref{eq:D1_telescopic_sum} leads to
	\begin{equation}
		\label{eq:D1_conclusion}
		J_{\infty}(x;\mu_{N,\hpara}, \tpara) \leq \frac{1}{1 - (\alpha^\ast_N(\epara) + \epsilon_N)}V_N(x; \hpara),
	\end{equation}
	which implies that $\ctx{\tpara}(\infty;x,\mu_{N,\hpara}) = 0$. Thus, the closed-loop system \eqref{eq:sec_prob:closed_loop_system} is asymptotically stable.
\end{proof}

\subsection{Proof of Theorem \ref{thm:affine_performance_ce_mpc}}
\label{appendix:D2}
\begin{proof}
	Noting \eqref{eq:parameter_perturbation} and $V_N(x; \tpara) \leq V_{\infty}(x; \tpara)$, \eqref{eq:D1_conclusion} implies
	\begin{equation}
		\label{eq:D2_bias_bound_1}
		J_{\infty}(x;\mu_{N,\hpara}, \tpara) \leq \frac{V_{\infty}(x; \tpara) + \beta_N(\epara; x)}{1 - (\alpha^\ast_N(\epara) + \epsilon_N)}.
	\end{equation}
	On the other hand, using \eqref{eq:sec_prob:mpc_linear_bound}, it holds that
	\begin{equation}
		\label{eq:D2_bias_bound_2}
		J_{\infty}(x;\mu_{N,\hpara}, \tpara) \leq \frac{V_{\infty}(x; \tpara) + \umpc\cost^\star(x)}{1 - (\alpha^\ast_N(\epara) + \epsilon_N)}.
	\end{equation}
	Combining \eqref{eq:D2_bias_bound_1} and \eqref{eq:D2_bias_bound_2} leads to \eqref{eq:performance_ce_mpc_affine}. In addition, if Assumption \ref{ass:local_scalable_input_perturbation_bound} holds, then by \eqref{eq:parameter_perturbation_linear} and $V_N(x; \tpara) \leq V_{\infty}(x; \tpara)$, \eqref{eq:D1_conclusion} implies that
	\begin{equation}
		\label{eq:D2_ratio_bound_1}
		J_{\infty}(x;\mu_{N,\hpara}, \tpara) \leq \frac{\big(1 + \beta^{\ast}_{N}(\epara)\big)V_{\infty}(x; \tpara)}{1 - (\alpha^\ast_N(\epara) + \epsilon_N)}.
	\end{equation}
	In addition, using \eqref{eq:sec_prob:mpc_linear_bound}, it also holds that
	\begin{equation}
		\label{eq:D2_ratio_bound_2}
		J_{\infty}(x;\mu_{N,\hpara}, \tpara) \leq \frac{\big(1 + \umpc \big)V_{\infty}(x; \tpara)}{1 - (\alpha^\ast_N(\epara) + \epsilon_N)}.
	\end{equation}
	Finally, \eqref{eq:performance_ce_mpc_linear} can be obtained by leveraging \eqref{eq:D2_ratio_bound_1} and \eqref{eq:D2_ratio_bound_2}.
\end{proof}

\section{Details of Section \ref{sec:examples}}
\label{appendix:E}

\subsection{More on Corollary \ref{coro:linear_quadratic_control}}
\label{appendix:E1}
Define $\Phi_{A}$ and $G_{A,B}$ as
\begin{equation*}
	\Phi_{A} := 
	\begin{bmatrix}
		\mathbf{I} \\
		A \\
		A^2 \\
		\vdots \\
		A^N
	\end{bmatrix},\quad
	G_{A,B} :=
	\begin{bmatrix}
		0 & 0 & \cdots & 0 \\
		B & 0 & \cdots & 0 \\
		AB & B & \cdots & 0 \\
		\vdots & \vdots & \vdots & \vdots \\
		A^{N-1}B & A^{N-2}B & \cdots & B
	\end{bmatrix},
\end{equation*}
and denote $\bar{Q}_{N} = \mathbf{I}_{N+1} \otimes Q$ and $\bar{R}_N = \mathbf{I}_{N} \otimes R$, where $\mathbf{I}_{N+1}$ is the identity matrix of dimension $N$ and $\otimes$ stands for the Kronecker product. The optimal input $\bar{u}^\star_N(x) = [[u^\star_{0|t}(x)]^\top, [u^\star_{1|t}(x)]^\top, \dots, [u^\star_{N-1|t}(x)]^\top]^\top$ for the unconstrained linear quadratic control problem is given by
\begin{equation*}
	\bar{u}^\star_N(x) = -(\bar{R}_N + G^\top_{A,B}\bar{Q}_{N}G_{A,B})^{-1}G^\top_{A,B}\bar{Q}_{N}\Phi_A x,
\end{equation*}
from which the inputs $u^\star_{k|t}(x) = K_{k,N}x, \;(k = 0,1,\dots,N-1)$ are obtained accordingly. It is clear that $K_{k,N}$ is a function of $\para$ through its dependency on $A$ and $B$. As established in \cite[Lemma 2]{mania2019certainty} and \cite[Lemma 14]{shi2025suboptimality}, there exists $L_{K} > 0$ such that $\|K_{k,N}(\para) - K_{k,N}(\hpara)\| \leq L_{K}\mpara$ for all $k \in \bI_{0:N-1}$. Moreover, without loss of generality, the local region $\Omega$ is chosen as $\mathcal{B}(0, r_{\mathrm{LQ}})$, where $r_{\mathrm{LQ}} = [\epsilon^\star_K/\emin{Q}]^{\frac{1}{2}}$ with $\epsilon^\star_K := \min_{\para \in \setpara, k \in \bI_{0:N-1}} \min_{i}\|[E_uK_{k,N}(\para)]^{\top}_{(i,:)}\|^{-2}_{Q^{-1}}$. As such, we define $\omelq := \min\left\{L_{K}, \frac{2}{(1-\edsr)r_{\mathrm{LQ}}}\right\}$, and $\sueds{k}^\ast$ in \eqref{eq:subsec_para_perturbation:refined_eds_separable} admits a linear upper bound $\sueds{k}^\ast(\mpara) \leq \omelq \mpara$ for all $k \in \bI_{0:N-1}$. The constant $\omelq$ is used to characterize the input perturbation in a compact way. Having established that both $\sueds{k}^\ast$ and $\ldd$ are linear functions of $\mpara$, the quadratic form (i.e., $f(z) = az^2 + bz$) of $\alpha^\ast_N$ and $\beta^\ast_N$ follows directly based on Propositions \ref{prop:parametric_state_perturbation} and \ref{prop:linear_parameter_perturbation}.

}

\section*{References}
\vspace*{-0.5cm}
{\footnotesize
\bibliographystyle{IEEEtran}
\bibliography{4-references/additional_reference,4-references/basic_references,4-references/books,4-references/learning,4-references/network,4-references/oco,4-references/intro_used,4-references/ourpaper}
}
\vspace*{-1cm}
\begin{IEEEbiography}[{\includegraphics[width=1in,height=1.25in,clip,keepaspectratio]{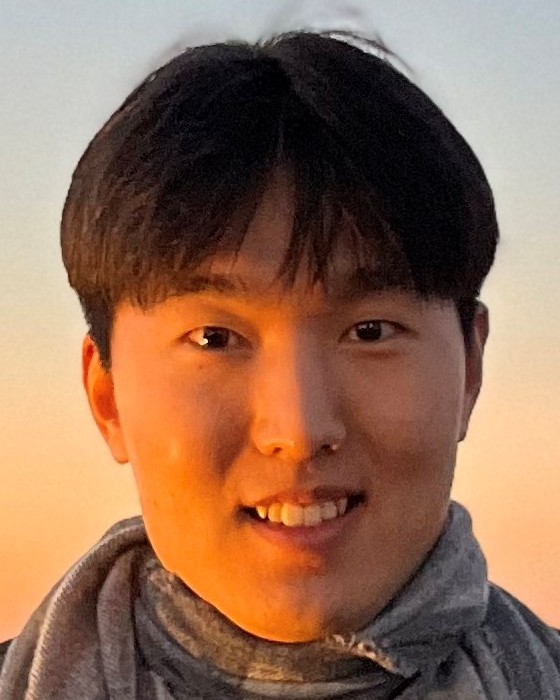}}]{Changrui Liu} received the BSc degree (cum laude) from Harbin Institute of Technology (HIT) in Harbin, China, he graduated from Honors School, HIT and was awarded Best 10 Graduate, 2019. He then received the MSc degree from TU Delft, Aerospace Engineering. He is currently pursuing the PhD degree at the Delft Center for Systems and Control, TU Delft. His research interests are performance guarantees in machine learning and optimization.
\end{IEEEbiography}
\vspace*{-1cm}
\begin{IEEEbiography}[{\includegraphics[width=1in,height=1.25in,clip,keepaspectratio]{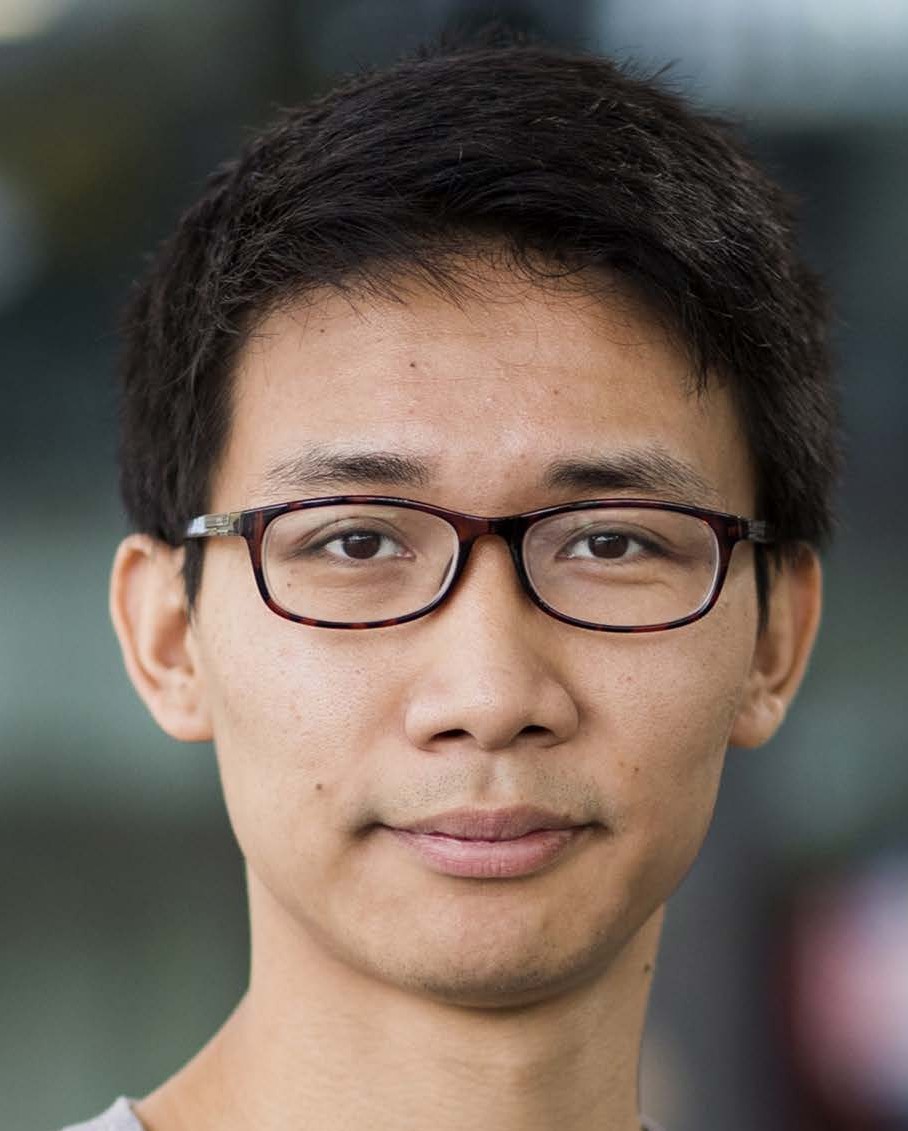}}]{Shengling Shi} is an assistant professor at Delft Center for Systems and Control, TU Delft, and he was a postdoctoral associate at the Massachusetts Institute of Technology, USA. He was also a postdoc at the Delft Center for Systems and Control, TU Delft, the Netherlands. He received the Ph.D. degree from the Eindhoven University of Technology in 2021. His research interests include system identification, model predictive control, and their applications.
\end{IEEEbiography}
\vspace*{-1cm}
\begin{IEEEbiography}[{\includegraphics[width=1in,height=1.25in,clip,keepaspectratio]{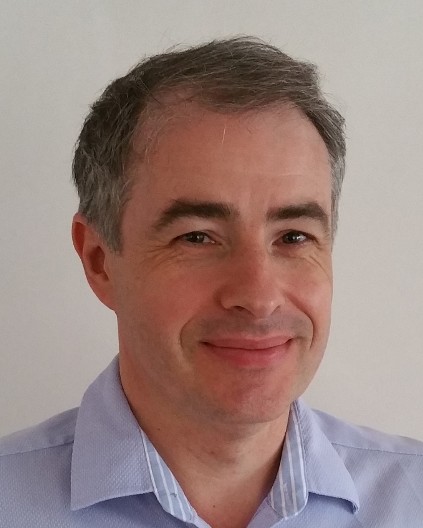}}]{Bart De Schutter} (IEEE member since 2008, senior member since 2010, and fellow since 2019) is a full professor and head of 
department at the Delft Center for Systems and Control of Delft 
University of Technology in Delft, The Netherlands. He is also an IFAC fellow.

Bart De Schutter is a senior editor of the IEEE Transactions on 
Intelligent Transportation Systems and was an associate editor of IEEE Transactions on Automatic Control.  His current research interests 
include integrated learning- and optimization-based control, multi-level and multi-agent control, 
and control of hybrid systems.
\end{IEEEbiography}

\end{document}